\documentclass{amsart}
\usepackage{amscd}
\numberwithin{equation}{section}
\let\cal\mathcal

\def\r{\rightarrow}

\newtheorem{lemma}{Lemma}[section]
\newtheorem{proposition}[lemma]{Proposition}
\newtheorem{theorem}[lemma]{Theorem}
\newtheorem{maintheorem}[lemma]{Main Theorem}

\newtheorem{question}[lemma]{Question}
\newtheorem{corollary}[lemma]{Corollary}

\newtheorem{fact}[lemma]{Fact}

\theoremstyle{definition}

\newtheorem{example}[lemma]{Example}
\newtheorem{definition}[lemma]{Definition}
\newtheorem{convention}[lemma]{Convention}
\newtheorem{mainconjecture}[lemma]{Conjecture 1}

\newtheorem{notation}[lemma]{Notation}

\newtheorem{conjectures}[lemma]{Conjectures}
\newtheorem{strongconjecture}[lemma]{Strong Conjecture}

{

\newtheorem*{prooftheoremA}{Proof of theorem A}

\newtheorem*{prooftheorem}{Proof of the theorem}

\newtheorem*{sketchproof}{Sketch of the proof}
\newtheorem*{proofproposition}{Proof of the proposition}
}

\theoremstyle{remark}

\newtheorem{remark}[lemma]{Remark}

\newcommand{\Acal}{\mbox{$\cal A$}}

\newcommand{\Gcal}{\mbox{$\cal G$}}

\newcommand{\Lcal}{\mbox{$\cal L$}}

\newcommand{\Zcal}{\mbox{$\cal Z$}}

\newcommand{\Ocal}{\mbox{$\cal O$}}
\newcommand{\Ucal}{\mbox{$\cal U$}}

\newcommand{\Wcal}{\mbox{$\cal W$}}
\newcommand{\Rcal}{\mbox{$\cal R$}}

\newcommand{\Scal}{\mbox{$\cal S$}}

\newdimen\uboxsep \uboxsep=1ex
\def\uboxn#1{\vtop to 0pt{\hrule height 0pt depth 0pt\vskip\uboxsep
\hbox to 0pt{\hss #1\hss}\vss}}

\def\uboxs#1{\vbox to 0pt{\vss\hbox to 0pt{\hss #1\hss}
\vskip\uboxsep\hrule height 0pt depth 0pt}}

\title[Set-theoretic solutions of the Yang-Baxter equation]{A combinatorial approach to the
set-theoretic solutions of the Yang-Baxter equation}
\thanks{The author was partially supported by the
Department of Mathematics of  Harvard University,
  by Grant MM1106/2001 of the Bulgarian National Science Fund
  of the Ministry of Education and Science,
  and by the Programme "Noncommutative Geometry" of the European Science Foundation}
\author{Tatiana Gateva-Ivanova}
\address{Institute of Mathematics and Informatics\\
Bulgarian Academy of Sciences\\
Sofia 1113, Bulgaria\\
and\\
American University in Bulgaria\\
2700 Blagoevgrad, Bulgaria
}

\email{ tatyana@aubg.bg, tatiana@math.bas.bg, tatiana@math.harvard.edu}

\begin{document}
\begin{abstract}
A  bijective map $r: X^2 \longrightarrow X^2$, where $X = \{x_1,
\cdots , x_n \}$ is a finite set, is called a \emph{set-theoretic
solution of the Yang-Baxter equation} (YBE) if  the braid relation
$r_{12}r_{23}r_{12} = r_{23}r_{12}r_{23}$ holds in $X^3.$  A
non-degenerate involutive solution $(X,r)$ satisfying $r(xx)=xx$,
for all $x \in X$,  is called \emph{square-free solution}. There exist
close relations between the square-free  set-theoretic solutions
of YBE, the semigroups of I-type, the semigroups of skew
polynomial type, and the Bieberbach groups, as it was first shown in a
joint paper with Michel Van den Bergh.

In this paper we continue the study of square-free solutions
$(X,r)$ and the associated Yang-Baxter algebraic structures --- the
semigroup $S(X,r)$, the group $G(X,r)$ and the $k$- algebra $A(k,
X,r)$ over a field $k$, generated by $X$ and with quadratic
defining relations naturally arising and uniquely determined by
$r$. We study the properties of the associated Yang-Baxter
structures and prove a conjecture of the present author that the
three notions: a square-free solution of (set-theoretic) YBE, a
semigroup of I type, and a semigroup of skew-polynomial type, are
equivalent. This implies that the Yang-Baxter algebra $A(k, X,r)$
is Poincar\'{e}-Birkhoff-Witt type algebra, with respect to some
appropriate ordering of $X$. We conjecture that every square-free
solution of YBE is retractable, in the sense of Etingof-Schedler.
\end{abstract}
\maketitle
\section{Introduction}

The Yang-Baxter equation appeared in 1967 \cite{Yang} in
Statistical Mechanics and turned out to be one of the basic
equations in mathematical physics, and  more precisely for
introducing the theory of quantum groups. At present the study of
quantum groups, and, in particular, the solutions of the
Yang-Baxter equation attracts the attention of a broad circle of
scientists and mathematicians.

Let $V$ be a vector space over a field $k$. We recall that a  linear automorphism $R$
of $V\otimes V$ is  \emph{a solution of the Yang-Baxter equation}, if
the equality
\begin{equation}
\label{YBE}
 (R\otimes id_V)(id_V\otimes R)(R\otimes id_V) =
 (id_V\otimes R)(R\otimes id_V)(id_V\otimes R)
\end{equation}
holds in the authomorphism  group of $V\otimes V\otimes V.$
$R$ is a solution of the \emph{quantum Yang-Baxter equation} (QYBE)
if
\begin{equation}
\label{QYBE}
R^{12}R^{13}R^{23} = R^{23}R^{13}R^{12}
\end{equation}
where $R^{ij}$ means $R$ acting on the i-th and j-th component.

Finding all solutions of
the Yang-Baxter equation is a difficult task far from being resolved.
Nevertheless many solutions of these equations have been found during
the last 20 years and the related algebraic structures
(Hopf algebras) have been studied (for example see \cite{K}).
Most of these solutions were "deformations" of the identity solution.
In 1990 V. Drinfeld \cite{D}  posed the problem
of studying a class of solutions that are obtained in a different
way - the so called \emph{set-theoretic solutions}.
\begin{definition}
\label{d1}
Let $X$ be a nonempty set.
Let \
$r:$ $X\times X$ $\longrightarrow$ $X\times X$
be a bijection of the Cartesian product
$X\times X$ onto itself. The map $r$ is called
\emph{a set-theoretic solution of the Yang-Baxter equation},
if
\[
(r\times id_X)(id_X\times r)(r\times id_X)=
(id_X\times r) (r\times id_X)(id_X\times r).
\]

Each set-theoretic solution $r$ of the Yang-Baxter equation
induces an operator $R$  on $V\otimes V$ for the vector space $V$
spanned by $X$, which is, clearly, a solution of  \ref{YBE}.
Various works dealing with set-theoretic solutions appeared during
the last decade, cf. \cite{W}, \cite{H}, \cite{TM}, \cite{ESS},
\cite{E1}, \cite{S}, \cite{Lu}, \cite{O}, \cite{Rump}.

The purpose of this paper is first to present some recent
conjectures on the set-theoretic solutions of the Yang-Baxter
equation, and to give an account of the research in this area,
and,  second to continue the study of the general algebraic and
homological properties of the algebraic structures related to the
so called square-free solutions. Our approach is combinatorial. To
each solution $(X, r)$ we associate a  semigroup $S= S(X, r)$, a
group $G= G(X, r)$ (the group was also studied in  \cite{ESS}),
and a quadratic algebra over a field $k$, $A(k, X, r)\simeq kS$,
each of them with a set of $n$ generators $X$ and  with quadratic
defining relations $\Re(X, r)$ naturally arising and uniquely
determined by $r$. We study the "behaviour"  of these relations,
and use the obtained information for establishing structural and
homological properties of the associated algebraic objects.
This approach is natural, for usual linear solutions one has similar ideas for
instance  Manin's work \cite{Manin}. In the case of set-theoretic
solutions to YBE it was initiated in the joint paper with Michel
Van den Bergh \cite{TM}, and applied to the study of the close
relations between different mathematical objects such as
set-theoretic solutions of the Yang-Baxter equation, semigroups of
I-type (which appeared in the study of Sklyanin algebras) and the
semigroups $S_0$ associated with the  class of  skew-polynomial
rings with binomial relations, introduced and studied in \cite{T1}
and \cite{T2}. The semigroups $S_0$, called \emph{ semigroups of
skew-polynomial type} are standard finitely presented, more
precisely, they are defined  in terms of a finite number of
generators and quadratic square-free relations, which form  a
Groebner basis (or equivallently, the algebra $A=kS$ is a PBW
algebra) cf. \ref{defskew}. It is proven in \cite{TM} that each
skew-polynomial semigroup
 $S_0$
defines a nondegenerate set-theoretic solution
$r=r(S_0)$ of the Yang-Baxter equation.
In connection with this result the present author made the
conjecture that under the restriction that $X$ is finite and
``square-free'' i.e.  $r(x,x)= (x,x)$ for each $x \in X,$
all nondegenerate involutive
solutions can be obtained in this way, cf. \ref{mainconjecture}.

In this work we will not be in a position to develop specific physical applications
but already we can say that several of the structures we introduce are
highly relevant for physics. For example, the groups $G(X,r)$ act on each other to
form a matched pair of groups and are hence a natural source of quantum
groups of bicrossproduct type. More details are to appear in our sequel  \cite{TShahn} .
Bicrossproduct quantum groups themselves  are increasing importance in noncommutative
geometry as for example the Connes-Kreimer quantum groups associated to renormalisation,
the $\kappa$-Poincar\'e quantum groups related to deformed spacetime,
 and the original 'Planck-scale' quantum group; see \cite{Shahn} for this background.

\section{Basic notions and results}

\label{sec2} In this section we first recall  some basic notions,
definitions, and results, from  \cite{ESS}, and \cite{TM}. They
are related to both quantum group theory and noncommutative
algebra, so we recall them for convenience of readers with various
mathematical background. Next we formulate the main results of the
paper and a conjecture about set-theoretic solutions of YBE.

We fix a finite nonempty set $X$ with $n$ elements.
We shall often identify the sets $X\times X$
and $X^2,$ the set of all monomials of length two in
the free semigroup $\langle X\rangle.$

\begin{definition}
\label{def1} \cite{ESS} Let $r: X\times X \r  X\times X $ be a
bijective map, we shall  refer to it as $(X, r)$. The
components of $r$ are the maps $\Lcal: X\times X \rightarrow  X$
and $\Rcal: X\times X  \rightarrow  X$  defined by the equality
\[
r(x,y) = (\Lcal _x(y), \Rcal _y(x)).
\]
(i) $(X, r)$ is \emph{left nondegenerate} if for each $x$ the map
$\Lcal _x(y)$ is a bijective function of $y$; $(X,r)$ is
\emph{right nondegenerate} if for each $y$ the map $\Rcal _y(x)$
is a bijective function of $x$;  $(X,r)$ is \emph{nondegenerate}
if it is left and right nondegenerate.

(ii) $(X,r)$ is \emph{involutive} if
\begin{equation}
\label{I}
r^2 = id_{X\times X}
\end{equation}

(iii) $(X,r)$ is \emph{a braided set} if $r$
satisfies the braid relation:
\begin{equation}
\label{B}
r_1r_2r_1=r_2r_1r_2,
\end{equation}
where
$r_1= r\times id_X$ and $r_2=id_X\times r.$

(iv)$(X, r)$ is \emph{symmetric} if it is braided and involutive.

(v) If $(X, r)$ is a braided, involutive and nondegenerate set
 we shall call it simply \emph{a solution}.
\end{definition}
Clearly, every braided set presents a set-theoretic solution of
the Yang-Baxter equation. A general study of
nondegenerate symmetric sets was given in \cite{ESS}.

In \cite{TM} was found a special class of
of solutions, here we call them  \emph{square-free solutions}
(cf. \ref{def2}),
which are defined via
the semigroups with relations of skew-polynomial type.
These semigroups were
introduced and studied first in \cite{T1}. The study continued in
\cite{T2}, \cite{T3}, \cite{TM},
\cite{JO}, cf. also  \cite{TJO}.

\begin{definition}
\label{def2} A map  $r: X^2 \r  X^2 $  is \emph{square-free} if it
acts trivially on $diag(X^2),$ i.e. $r(xx) = xx,$ for all $x \in
X.$
\end{definition}

\begin{example}
\label{ex1}

Let $X$ be a nonempty set and let $r(xy) = yx.$ Then $(X, r)$
is a square-free solution, which is called
\emph{the trivial solution}.
\end{example}

\begin{example}(\emph{Permutational solution}, Lyubashenko, \cite{D}).
\label{ex2}
Let $X$ be a non-empty set, let $f, g$ be maps
$X \r X$ and let $r(xy) = g(y)f(x).$ Then a) $(X, r)$ is nondegenerate
if and only if $f$ and $g$ are bijective; b) $(X, r)$ is braided
if and only if $fg=gf$; c) $(X,r)$ is involutive if and only if
$f = g^{-1}$.
\end{example}

\begin{remark}
\label{rem1}
Note that for any permutation $f$ of $X$, the map
$r$ defined as $r(xy) = f(y)f^{-1}(x)$,
is a solution, but in general $r$
is not square-free. In fact, a permutational involutive solution $r$ is
square-free if and only if
$f= id_{X},$ i.e. $r= id_{X^2}$.
Nevertheless, we prove in \ref{CC}  that each square-free solution
behaves "locally"
as a permutational solution.
\end{remark}
Clearly, when the order  $\mid X\mid =2$, the only square-free
solution $(X, r)$ is the trivial one. The lowest order of $X$
which allows a nontrivial, square-free solution is $3$, as shows
the following.
\begin{example}
\label{ex3} Let $X= \{ x_1,x_2,x_3 \}$. Up to re-numerating of the
set $X$ there exists a uniqie non-trivial square-free solution
$(X, r)$ namely :
\[ r(x_3x_1) = x_2x_3, \ r(x_2x_3)= x_3x_1;
\]
\[
 r(x_3x_2)=x_1x_3, \ r(x_1x_3) = x_3x_2,
\]
\[
r(x_2x_1) = x_1x_2, \ r(x_1x_2)= x_2x_1, \ r(x_ix_i)=x_ix_i,
i=1,2,3.
\]
\end{example}
Up to isomorphism of solutions, there exist $5$ square-free
solutions $(X, r)$ with $\mid X \mid=4$. The one with the greatest
number nontrivial relations is given in the following example.
\begin{example}
\label{example4}
 Let $X= \{ x_1,x_2,x_3,x_4 \}$ and let $r$ be defined
as:
\[
r(x_1x_3) = x_4x_2, r(x_4x_2)= x_1x_3,
r(x_1x_4) = x_3x_2, r(x_3x_2)=x_1x_4,
\]
\[
r(x_2x_3)= x_4x_1, r(x_4x_1)=x_2x_3,
r(x_2x_4)=x_3x_1, r(x_3x_1)=x_2x_4,
\]
\[
r(x_1x_2)= x_2x_1,
r(x_2x_1)=x_1x_2,
r(x_3x_4)=x_4x_3,
r(x_4x_3)=x_3x_4,
\]
\[
r(x_ix_i) = x_ix_i, i = 1,\cdots, 4.
\]
Then $(X, r)$ is a square-free solution.
Consider the permutation $\sigma = (12)(34)$.
For $x, y$ which belong to different orbits of $\sigma$
one has $r(xy)=\sigma(y)\sigma^{-1}(x),$
and when $x$ and $y$ belong to the same orbit,
then $r(xy)=\sigma^2(y)\sigma^{-2}(x)=yx.$
\end{example}
\begin{definition}
\label{def3}
The braid group $B_n$ is the group generated by $n$ generators
$b_1, \cdots ,  b_n$ and defining relations
\begin{equation}
\label{Br}
b_ib_j = b_jb_i,    \mid i-j \mid > 1;
\end{equation}
\begin{equation}
\label{Br1}
b_ib_{i+1}b_i = b_{i+1}b_ib_{i+1}.
\end{equation}
\end{definition}
Recall that the symmetric group $S_n$ is isomorphic to the quotient
of $B_n$ by the relations $b_i^2 = 1.$

The following remark is obvious, see for example \cite{ESS}.
\begin{remark}
\label{P1}
Let $m \geq 3$ be an integer.
 (i) The assignment $b_i \r r^{ii+1}, 1\leq i \leq m-1,$ extends to an action of
$B_m$ on $X^m$  if and only if $(X, r)$ is a braided set. (ii) The
assignment $b_i \r r^{ii+1}, 1\leq i \leq m-1,$ extends to an
action of $S_m$ on $X^m$  if and only if $(X, r)$ is a symmetric
set. (Here, as usual, $r^{ii+1}= id_{X^{(i-1)}} \times r \times
id_{X^{(m-i-1)}}$).
\end{remark}
The next well-known fact (see \cite{ESS}) gives the relation
between the braided sets (i.e. the set-theoretic solutions of the
Yang-Baxter equation) and the set-theoretic solutions of the
quantum Yang-Baxter equation.
\begin{fact}
\label{P2} Let $r: X^2 \r  X^2 $ be a bijection, $\sigma: X^2 \r
X^2 $ be the flip $\sigma(xy) = yx$, for all $x, y \in X.$ Let $R
= \sigma \circ r.$ (i.e. $R$ is the so called \emph{R-matrix}
corresponding to $r$). Then $r$ satisfies the set-theoretic
Yang-Baxter equation
 if and only if $R$ satisfies the quantum Yang-Baxter equation:
\begin{equation}
\label{QYBE1}
R^{12}R^{13}R^{23} = R^{23}R^{13}R^{12}.
\end{equation}
Furthermore, $r$ is involutive if and only if $R$ satisfies
\ref{QYBE1} and the unitarity condition
\begin{equation}
\label{eq2}
R^{21}R = 1.
\end{equation}
\end{fact}
In the spirit of a recent trend called a $\emph{combinatorial
approach in algebra},$ to each bijective map $r: X^2 \r  X^2$ we
associate canonically finitely presented algebraic objects (see
precise definition in \ref{associatedobjects}) generated by $X$
and with quadratic defining relations $\Re$ naturally determined
as
\begin{equation}
\label{defrelations}
 \Re=\Re(r) =\{(u=r(u))  \mid u \in X^2, u \neq
r(u) \; \text{as words in} X^2\}
\end{equation}

We study the close relations between the combinatorial properties
of the defining relations, e.g. of the map $r$, and the structural
 properties of the associated algebraic objects.
\begin{notation}
\label{freealgebra}
For a non-empty set $X$, as usual, we denote by
$\langle X \rangle$ the free semigroup generated by $X,$
and by $k\langle X \rangle$-
 the free associative $k$-algebra
generated by $X$, where $k$ is an arbitrary field.
For a set $F \subseteq k\langle X \rangle$, $(F)$ denotes
the two sided ideal of $k\langle X \rangle$, generated by $F$.
\end{notation}

\begin{definition}
\label{associatedobjects} Assume that $r:X^2 \longrightarrow X^2$
is an involutive, bijective map.

(i) The semigroup
\[
S =S(X, r) = \langle X; \Re(r) \rangle,
\]
 with a
set of generators $X$ and a set of defining relations $ \Re(r),$
is called \emph{the semigroup associated with $(X, r)$}.

(ii) The \emph{the group $G=G(X, r)$ associated with} $(X, r)$ is
defined as
\[
G=G(X, r)={}_{gr} \langle X; \Re (r) \rangle.
\]

(iii) For arbitrary fixed field $k$, \emph{the
$k$-algebra associated with} $(X ,r)$ is defined as
\begin{equation}
\label{Adef} \Acal = \Acal(k,X,r) = k\langle X \rangle/(\Re(r)).
\end{equation}
\end{definition}
Clearly $\Acal$ is a quadratic algebra, generated by $X$ and
 with defining relations  $\Re(r).$
Furthermore, $\Acal$  is isomorphic to the semigroup algebra $kS(X, r).$

Manin, \cite{Manin}, introduced the notion of
a \emph{Yang-Baxter algebra}.
He calls
\emph{a Yang-Baxter algebra}
a quadratic  algebra  $A$ with defining
relation determined via arbitrary  fixed Yang-Baxter operator.
In this spirit we give the following definition.
\begin{definition}
Assume $(X, r)$ is a solution.
Then  $S(X,r),$ $G(X, r)$ and $\Acal(k,X,r)$
are called respectively: \emph{the  Yang-Baxter
semigroup}, \emph{the  Yang-Baxter group},
 and
\emph{the  Yang-Baxter $k$-algebra, associated to} $(X,r)$.
We shall also use the abbreviation "\emph{YB}" for
"Yang-Baxter".
\end{definition}
In
the case when $(X, r)$ is a solution, $G(X, r)$ is also called
 the \emph{the structure group of $(X, r)$}, see \cite{ESS}.
\begin{example} Let $(X, r)$ be the trivial solution, i. e. $r(xy)=yx$,
for all $x, y \in X,$ then clearly, $S(X, r)=[x_1, \cdots, x_n]$,
is the free abelian semigroup generated by $X$, $G(X, r)=Z^X$, is
the free abelian group generated by $X$, and $\Acal(k,X,r)= k[x_1,
\cdots, x_n]$ is the commutative polynomial ring over $k$.
\end{example}
\begin{definition}
\label{YBS} Let $S= \langle X; \Re \rangle $ be a semigroup with a
set of generators $X$ and a set of quadratic binomial defining
relations:
\[
\Re=\{ xy=y^{\prime}x^{\prime} \mid x, y, x^{\prime}, y^{\prime}
\in X \},
\]
We assume that each monomial $u \in X^2$, occurs in at most one
relation in $\Re.$ Define the map $r=r(S): X^2 \rightarrow  X^2$
as follows:

(i) $r(xy) = xy $, if $xy$ is a monomial of
length 2 which does not occur
 in any relation in $\Re$;
and

(ii) if $(xy=y^{\prime}x^{\prime}) \in \Re$, then we set
$r(xy)=y^{\prime}x^{\prime}$ and $r(y^{\prime}x^{\prime})=xy.$

We call $r(S)$ \emph{the map associated with the semigroup $S$}.
\end{definition}
Note that if $r$ is the map defined by the set of relations of a
YB- semigroup $S=\langle X; \Re \rangle$, then the set $(X; r)$ is
always symmetric, since clearly, $r^2=id_{X^2}$.

We give now an example of a Yang-Baxter semigroup $S$ with 11
generators. In fact, $S$ belongs to the class of semigroups of
skew-polynomial type, \ref{defskew}, and the  map $r(S)$ is a
square-free solution.
\begin{example}
\label{ex4} Let $S= \langle X; \Re \rangle$, where the set of
generators is $X= \{ 1, 2, \cdots , 8, a,b,c \}$ and  the defining
relations are:
\[
1a=a2, 2a=a1, 2b=b3, 3b=b2,  3a=a4, 4a=a3, 4c=c1, 1c=c4,
\]
\[
5a=a6, 6a=a5, 6b=b7, 7b=b6,  7a=a8, 8a=a7, 8c=c5, 5c=c8,
\]
\[
1b=b5, 5b=b1, 2c=c6, 6c=c2,  3c=c7, 7c=c3, 4b=b8, 8b=b4,
\]
\[
ab=ca, ac=ba, bc=cb, ij=ji, 1 \leq i,j \leq 8.
\]
\end{example}
\begin{remark}
Let $S_0$ be a semigroup of skew-polynomial type (see
\ref{defskew}). Let $r=r(S_0)$ be the map defined by the relations
of $S_0$. Then $(X, r)$ is a square-free solution (cf \cite{TM},
Th. 1.2,  also Theorem \ref{theoremA}). Furthermore, $S_0$ is a
cancellative semigroup, and has a group of quotients $gr(S_0)$,
which is a central localization of $S_0$, see \cite{JO}. It is
clear, that the groups $gr(S_0)$ and the associated  group $G(X,
r)$  are isomorphic. Moreover, the set $X$ is embedded in $G(X,
r).$
\end{remark}

The  semigroups of skew-polynomial type were discovered while the
author was searching for a new class of Artin-Schelter regular
rings. It turned out that the skew-polynomial rings with binomial
relations introduced and studied in \cite{T1}, \cite{T2},
\cite{T3} provide a class of Artin-Schelter regular rings of
arbitrary global dimension. Furthermore, with each ring $\Acal_0$
of this type we associate (uniquely) a semigroup $\Scal _0$ which
defines (via its relations) a non-degenerate set-theoretic
solution $r(\Scal _0)$ of the Yang-Baxter equation, cf. \cite{TM}.
It is easy to generalize this result by showing that each
skew-polynomial ring with binomial relations defines a solution of
the classical Yang-Baxter equation, see Theorem
\ref{skewpolyrings}. The semigroup $\Scal _0$ is called a
semigroup of skew-polynomial type. The results in \cite{TM} and
further study of the combinatorial properties of the solutions
inspired the following Conjecture, which we reported first in a
talk at the International Conference in Ring Theory, Miskolc 1996,
see also \cite{TConstanta}, \cite{TLovetch}.
\begin{mainconjecture} \cite{T5}
\label{mainconjecture} Let $(X, r)$ be a square-free (non-degenerate,
involutive) solution of the Yang-Baxter equation. Then the set $X$
can be ordered so, that the associated semigroup $S= S(X, r)$ is
of skew-polynomial type.
\end{mainconjecture}
\begin{definition}
\label{defskew} We say that the semigroup $S_0$ is \emph{a
semigroup of skew-polynomial type}, (or shortly, \emph{a
skew-polynomial semigroup}) if it has a standard finite
presentation as $S_0 = \langle X; \Re_0 \rangle $, where the set
of generators $X$ is ordered: $x_1 < x_2 < \cdots < x_n,$ and the
set
\[
\Re _0 =\{ x_jx_i=x_{i^{\prime}}x_{j^{\prime}})\mid 1 \leq i< j
\leq n, 1 \leq i^{\prime} < j^{\prime} \leq n \},
\]
contains precisely $n(n-1)/2$ quadratic square-free binomial
defining relations, each of them satisfying the following conditions:

i) each monomial $xy\in X^2$, with $x\neq y$, occurs in exactly
one relation in $\Re _0$; a monomial of the type $xx$ does not
occur in any relation in $\Rcal _0$;

ii) if  $(x_jx_i=x_{i^{\prime}}x_{j^{\prime}})\in \Rcal _0$, with
$1 \leq i< j  \leq n,$ then  $i^{\prime} < j^{\prime}$, and $j >
i^{\prime}$.

[ further studies show that this also implies $i < j^{\prime}$
 see \cite{T2}]

iii) the monomials $x_kx_jx_i$ with $k>j>i, 1\leq i,j,k, \leq n$
do not give rise to new relations in $S_0$, or equivalently,  cf.
\cite{B}, $\Re _0$ is a Groebner basis with respect to the
degree-lexicographic ordering of the free semigroup $\langle X
\rangle$.
\end{definition}
\begin{remark}
\label{PBW} Suppose $S_0$ is a semigroup of skew-polynomial type.
It follows from the Diamond Lemma \cite{B} that, each element $w$
of $S$ can be presented uniquely as an ordered monomial
\[ w =x_1^{\alpha_1}x_2^{\alpha_2}\cdots x_n^{\alpha_n}\]
where $ \alpha_i \geq 0, 1\leq i \leq n$. This presentation is
called \emph{the normal form of} $w$ and denoted as $Nor(w).$ It
follows from the Diamond lemma, that two monomials $w_1, w_2$ in
the free  semigroup $\langle X \rangle$ are equal in $S$ if and
only if their normal forms coincide, $Nor(w_1)=Nor(w_2).$
Thus $S_0$ can be identified as a set with the set of
ordered monomials
\begin{equation}
\label{orderedmonomials}
\mathcal{N}_0=\{
x_1^{\alpha_1}x_2^{\alpha_2}\cdots x_n^{\alpha_n}\mid \alpha_i
\geq 0, 1\leq i \leq n\}.
\end{equation}
Furthermore, for an arbitrary field  $k$, the set $\mathcal{N}_0$
is a $k$- basis of the quadratic algebra
\[A_0 = k\langle X \rangle/(\Re_0 ) \simeq kS_0.
\] Clearly, $A_0$ is a
Poincar\'{e}-Birkghoff-Witt - algebra in the sense of Priddy
\cite{priddy} with $\mathcal{N}_0$ as a PBW-basis.
\end{remark}
\begin{remark}
 In \cite{JO}  the skew-polynomial semigroups $S_0$
 are called \emph{binomial semigroups}.
\end{remark}
We now recall the definition of the semigroups of $I$-type, see
\cite{TM}, which are closely related to both- the semigroups of
skew-polynomial type and the set-theoretic solutions of
Yang-Baxter equation. The rings of  I-type were introduced and
studied by J.Tate, and  M. Van den Bergh in their work on the
homological properties of Sklyanin Algebras, \cite{Michel-Tate}.
\begin{notation}
\label{Ucal} Till the end of the paper we shall denote by
\begin{equation}
\Ucal=[u_1, \cdots , u_n],
\end{equation}
the free commutative multiplicative semigroup generated by $u_1,
\cdots , u_n$.
\end{notation}
\begin{definition}
\label{Itypedef} \cite{TM}, A semigroup $S$ generated by $\{x_1,
\cdots , x_n\}$ is said to be of \emph{(left) I-type} if there
exists a bijection $ v: \Ucal \longrightarrow S$ called (\emph{a
left
  I-structure}),
such that $v(1)=1,$ and such that for each $a\in \Ucal$ there is
an equality of sets $\{v(u_1a), v(u_2a), \cdots , v(u_na)\}=
\{x_1v(a),x_2v(a), \cdots , x_nv(a)\}$. Analogously one defines
\emph{a right I-structure} $v_1:\Ucal \longrightarrow S$.
\end{definition}
\begin{remark} It can be extracted from \cite{TM},
see also \ref{Itype1}, that if $(X, r)$ is a square-free solution,
and $S= S(X,r)$ the associated YB semigroup, then

a) There exists a unique left $I$-structure
 $v: \Ucal \rightarrow S,$ \
  such that $v(u_i)=x_i,$ \ for $1\leq i \leq n.$

b) There exists a unique right $I$-structure
 $v_1: \Ucal \rightarrow S,$
  such that $v_1(u_i)=x_i,$ \ for $1\leq i \leq n.$
\end{remark}

In section \ref{latticestructure}, Proposition \ref{lattice}, we
show that a semigroup of $I$-type is a distributive lattice with
respect to the order induced from "one-sided" divisibility,
defined bellow.
\begin{definition}
\label{ldivisibilitydef}
For every pair $a, b \in S$ we set:

(i) $ a \mid_l b$, if and only if there exists a monomial
 $c\in S$, such that  $b = ca.$ We call this relation
 \emph{divisibility with respect to the left
multiplication}.

(ii) $ a \mid_r b$, if and only if there exists a monomial
 $c\in S$, such that  $b = ac.$ This relation is called
 \emph{divisibility with respect to the right
multiplication}.
\end{definition}


The following theorem proved in section \ref{equivalentnotions}
verifies Conjecture \ref{mainconjecture}.
\begin{maintheorem}
\label{theoremA} Assume that $X$ is a finite set of order $n \geq 1,$ and
$r:$ $X\times X \longrightarrow$ $X\times X$ is a square-free involutive
bijection.  Let $S =S(X, r)$ be the semigroup associated with $(X, r),$ and let
$\Acal=\Acal(k,X,r)$ be the quadratic $k$-algebra associated with
$(X, r)$, where $k$ is an arbitrary field. Then the following
conditions are equivalent.
\begin{enumerate}
\item \label{theoremA1}$(X, r)$ is  non-degenerate
solution of the set-theoretic Yang-Baxter equation.
 \item
 \label{theoremA2}
 $S =S(X, r)$ is a
semigroup of $I$-type.
 \item
 \label{theoremA3}
 There exists an ordering on $X,$
$X=\{x_1 < x_2 < \cdots < x_n \}$, such that $S= S(X, r)$ is a
semigroup of skew-polynomial type. \item \label{theoremA4} There
exists an ordering on $X,$ $X= \{x_1 < x_2 < \cdots < x_n \}$ such
that for every field $k$ the quadratic $k$-algebra $\Acal=\Acal
(k, X,r)$ is a Poincar\'{e}-Birkhoff-Witt algebra, with a
$k$-basis - the set of ordered monomials $\mathcal{N} _0.$
\end{enumerate}
Moreover, each of these conditions implies that the solution $(X,
r)$ is decomposable, i.e. $X$ a disjoint union of two nonempty
$r$-invariant subsets.
\end{maintheorem}
\begin{corollary}
Let $(X, r)$ be a square-free solution, with associated semigroup
$S=S(X, r)$. Then $(S, \mid_l)$
 is a distributive lattice.
Furthermore the left $I$-structure $ v: \Ucal \longrightarrow S$
is an isomorphism of lattices.
\end{corollary}

Condition \ref{theoremA}.\ref{theoremA2} implies cf. \cite{TM},
various nice algebraic an homological properties of the algebra
$\Acal=\Acal(k, X, r)$,  like being a Noetherian domain, Koszul,
Cohen-Macaulay, Artin-Schelter regular, etc. In particular the
semigroup $S$ is cancellative. Hence it is naturally embedded in
its group of quotients $gr(S)= G(X, r)$. We recall these results
in Theorem \ref{general}.

My student, M.S. Garcia Roman has shown that for an explicitly
given solution $(X, r),$ condition \ref{theoremA}.\ref{theoremA3}
is equivallent to a standard problem from Linear Programming.

In  \cite{TShahn} is presented a matched pairs approach to the
set-theoretic solutions of the Yang-Baxter equation. One of the
main results in \cite{TShahn}, given here as Theorem
\ref{theoremB} covers all known constructions of solutions $(X,
r)$, restricted to the case of square-free solutions, with $X$ a
finite set.



In section \ref{multipermutation} we study the generalized twisted
unions of solutions, and multipermutation solutions.

Section \ref{binomialYBE} gives an application of the Main Theorem to
a particular class of solutions of the classical Yang-Baxter equation.

 We close this section with the following conjecture
\begin{strongconjecture}
I. Every square-free solution $(X, r)$, where $X$ is a finite set
of order $n \geq 2$,  is retractible. Furthermore  $(X, r)$ is a
multipermutation solution of level $m < n.$

II. Every multipermutation square-free solution of level $m$ is a
generalized twisted union of multipermutation solutions of levels
$\leq m-1$
\end{strongconjecture}

\section{The cyclic condition and combinatorics in $S(X, r)$}
\label{cycliccondition}

In this section we introduce a combinatorial technique for
 non-degenerate square-free solutions $(X,r)$, which associates
cycles in $Sym(X)$ to each pair of elements $y, x$ in $X$. We call
the corresponding property of $r$ \emph{cyclic condition}. The
cyclic condition is the base for all combinatorial techniques in
this paper. We use it here to deduce more precise pictures of the
left and right actions of  the group $G(X,r)$ on $X$, and  to to
show that each involutive square-free solution acts ``locally'' as
a permutational solution. We obtain some important relations of
higher degrees in $S(X, r)$, and use the lengths of the cycles
occurring in $S(X, r)$ to associate an invariant integer $M=M(X, r)$
with every solution $(X,r)$
called
\emph{the cyclic degree} of $(X,r)$.

\begin{definition}
\label{CCdef} Let $r: X\times X \longrightarrow X\times X$ be a
bijection.
\begin{enumerate}
\item \label{weakCC} We say that $(X, r)$ satisfies \emph{the weak
cyclic condition}, if  for every pair $y,x \in X$, there exist two
disjoint cycles $\Lcal _{y}^x= (x_1, \cdots , x_m)$ and
$\Rcal_{x}^y= (y_k, \cdots , y_1)$ in the symmetric group
$Sym(X)$, such that $x = x_1, y=y_1,$, and for all $1 \leq i \leq
m, \ \ 1\leq j \leq k,$ there are equalities:
\begin{equation}
\label{rji} r(y_jx_i)= \Lcal
_{y}^x(x_{i})\Rcal_{x}^y(y_{j})=x_{i+1}y_{j-1},
\end{equation}
where $x_{m+1} := x_1$, and $y_0 := y_k$.

In particular, $r(y x)= \Lcal _{y}^x(x)\Rcal_{x}^y(y)=x_2y_k.$
\item \label{strongCC} $(X,r)$ satisfies \emph{the cyclic
condition}, if for every pair $y,x \in X$, there exist two
disjoint cycles $\Lcal _y^x= (x_1, \cdots , x_m)$ and $\Lcal
_x^y=(y_1, \cdots , y_k)$ in $Sym(X)$, such that $x = x_1,
y=y_1,$, and for all $1 \leq i \leq m, \ \ 1\leq j \leq k,$ there
are equalities:
\begin{equation}
\label{rij} r(x_iy_j)=y_{j+1}x_{i-1} \  \text{and}\ \  r(y_jx_i)=
x_{i+1}y_{j-1},
\end{equation}
where $x_0=x_m, x_{m+1}:= x_1, $, and $y_{0}:= y_k, y_{k+1}=y_1$.

In particular, for every pair $(y,x) \in X\times X$, the disjoint
cycles $\Lcal _y^x$ and $\Lcal _x^y$ satisfy:
\begin{equation}
r(y,x)= \Lcal _y^x (x)(\Lcal _x^y) ^{-1}(y), \ \text{and} \ \
r(x,y)= \Lcal _x^y (y)(\Lcal _y^x) ^{-1}(x).
\end{equation}
\end{enumerate}

We call $\Lcal _x^y$ and $\Lcal _y^x$ \emph{the pair of cycles
associated to} $(y, x).$
\end{definition}
\begin{remark}
Clearly,  the (strong) cyclic condition implies that $r$ is
involutive. We will show that every involutive square-free
solution $(X, r)$ satisfies the cyclic condition and use this to
study the left (and the right) action of $G(X,r)$ on $X$.  Note
that if the cyclic condition holds, and we set
\[
\sigma = \sigma_{y, x}= \sigma_{x, y}=(x_1, \cdots , x_m)(y_1,
\cdots , y_k) \in Sym(X),
\]
 the map $r$ is expressible "locally" as
a permutational solution
\[
r(y_jx_i) = \sigma (x_i)\sigma ^{-1}(y_j) \ \ \text{and} \  \
r(x_iy_j) = \sigma (y_{j})\sigma ^{-1}(x_i).
\]
If we do not assume involutiveness  for $r$, then,  in general,
only the weak cyclic condition is satisfied. We give an example,
see \ref{exnoninvolutive}, of a non-involutive solution in which
the cyclic condition does not hold.
\end{remark}
\begin{example}
\label{exnoninvolutive} Let $X= \{x_1, x_2, x_3, x_4, x_5, x_6 \}
$, and suppose the map $r: X^2 \rightarrow  X^2$ is defined as:
\[
x_1x_2 \leftrightarrow x_2,x_1 ; x_3x_4 \leftrightarrow  x_4x_3;
\]
\[
x_3x_5 \leftrightarrow x_5x_3 ; x_3x_6 \leftrightarrow x_6x_3;
\]
\[
x_4x_5 \leftrightarrow x_5x_4 ; x_4x_6 \leftrightarrow x_6x_4; xx
\leftrightarrow xx, \  \text{for all} \  x \in X
\]
\[
x_1x_3 \rightarrow  x_4x_2 \rightarrow x_1x_5 \rightarrow x_6x_2
\rightarrow x_1x_3 ;
\]
\[
x_1x_4 \rightarrow  x_3x_2 \rightarrow x_1x_6 \rightarrow x_5x_2
\rightarrow x_1x_4 ;
\]
\[
x_2x_3 \rightarrow  x_4x_1 \rightarrow x_2x_5 \rightarrow x_6x_1
\rightarrow x_2x_3 ;
\]
\[
x_2x_4 \rightarrow  x_3x_1 \rightarrow x_2x_6 \rightarrow x_5x_1
\rightarrow x_2x_4.
\]
Then $(X, r)$ is a non-involutive solution, with $r^4= id_{X^2}.$
Furthermore
\[
\Lcal_{x_1}= (x_3 x_4)(x_5x_6); \Rcal_{x_1} =(x_3x_6)(x_4x_5),  \
\text{and} \  \Rcal_{x_1} \neq (\Lcal_{x_1})^{-1}.
\]
\end{example}
Recall first a well known fact from \cite{ESS}.
\begin{fact}\cite{ESS}
\label{actions} Let $(X,r)$ be nondegenerate, $G=G(X, r)$. Then
$(X, r)$ is a braided set if and only if the following three
conditions are satisfied:
\begin{enumerate}
\item
The assignment $x \rightarrow \Lcal _x$ induces a left action of $G$ on $X$;
\item
The assignment $x \rightarrow \Rcal _x$ induces a right action of $G$ on $X$;
\item
The following equality holds for any $x, y, z \in X$:
\begin{equation}
\Lcal _{\Rcal _{\Lcal _y(z)}(x)}(\Rcal _z(y))=
\Rcal _{\Lcal _{\Rcal _y(x)}(z)}(\Lcal _x(y)).
\end{equation}
\end{enumerate}
\end{fact}
\begin{notation}
\label{Gorbits}  We shall denote by  $O_G(x)$ the orbit of $x,
x\in X,$ under the left action of $G$ on $X$ .
\end{notation}

\begin{lemma}
\label{lcc} With notation being as in \ref{CCdef},
\begin{enumerate}
\item $(X, r)$ satisfies the weak cyclic condition if and only if
for all $i,j, 1 \leq i \leq m$,  $1 \leq j \leq k$,  there are
equalities
\begin{equation}
\Lcal _{y_j}^{x_i}= \Lcal _y^x =(x_1, \cdots , x_m), \  \text{and}
\ \Rcal _{x_i}^{y_j}= \Rcal _x^y =(y_1, \cdots , y_k).
\end{equation}
 \item $(X, r)$ satisfies the cyclic
condition if and only if for all $i,j, 1 \leq i \leq m$, $1 \leq j
\leq k$ there are equalities
\begin{equation}
\label{cyclicconditiongeneral}
\Lcal _{y_j}^{x_i}= (\Rcal
_{y_j}^{x_i})^{-1} = (x_1, \cdots , x_m),
\end{equation}
and
\begin{equation}
 \Lcal _{x_i}^{y_j} = (\Rcal _{x_i}^{y_j})^{-1}=(y_1, \cdots , y_k).
\end{equation}
\end{enumerate}
\end{lemma}

The following theorem gives an account of various conditions on
the bijective maps
$r: X^2 \rightarrow X^2$ and the corresponding semigroup $S(X, r)$.
For some of them we assume neither that $r$ is
necessarily a solution of the Yang-Baxter equation,
nor we assume that $r$ is involutivene.

\begin{theorem}
\label{CC}
Let  $r: X^2 \rightarrow X^2$ be a bijective map, denoted by $(X, r)$.
Let $S=S(X, r)$ be the semigroup associated to $(X, r)$.
Let $\Lcal _x$ and $\Rcal _x$ be the left and right components of $r$,
introduced in \ref{def1}.
Consider the following conditions:
\begin{enumerate}
\item \label{L11} a) $(X, r)$ is left nondegenerate; b) $(X, r)$
is right nondegenerate. \item \label{ore} a) (Right Ore condition)
For every pair $a, b\in X$ there exists a unique pair $x,y\in X$,
such that $ax=by$; b) (Leftt Ore condition) For every pair $a,
b\in X$ there exists a  unique pair $z,t\in X$, such that $za=tb$.
\item \label{squarefree1} $(X,r)$ is square-free and
nondegenerate. \item \label{squarefree2} $\Lcal _x$ is a bijection
and $\Lcal _x(y) \neq x,$ for each $y\neq x$; $\Rcal _y$ is a
bijection and $\Rcal _y(x) \neq y,$ for each $y\neq x.$
\end{enumerate} Then the following is true:

A. The conditions \ref{L11} a),  and \ref{ore} a) are equivalent;
the conditions \ref{L11} b),  and \ref{ore} b) are equivalent;

B. The conditions \ref{squarefree1} and \ref{squarefree2} are equivalent.

C. If $(X,r)$ is a non-degenerate square-free solution of the
Yang-Baxter equation, (not necessarily involutive) then the weak
cyclic condition \ref{CCdef}\ref{weakCC} holds.

D. If $(X,r)$ is a non-degenerate involutive square-free solution
of the Yang-Baxter equation, then the cyclic condition
\ref{CCdef}.\ref{strongCC} holds.
\end{theorem}
\begin{proof}
A. $(\ref{L11}.a) \Longrightarrow (\ref{ore}.a))$ Let $a,b\in X$.
By our assumption the function $\Lcal _a$ is a bijection of $X$
onto itself, so there exists a unique $y$ such that $\Lcal
_a(y)=b$, hence the equality $r(ay)= \Lcal _a(y)\Rcal _y(a)$ gives
$r(ay)=bz,$ for some $z\in X$. But $r$ is a bijective map on $X^2$
onto itself, so $z$ is also determined uniquely. The implication
$(\ref{L11}.b) \Longrightarrow (\ref{ore}.b))$ is analogous.

The implications $(\ref{ore}.a) \Longrightarrow (\ref{L11}.a))$
and $(\ref{ore}.b) \Longrightarrow (\ref{L11}.b))$ are obvious.

B. $\ref{squarefree1} \Longrightarrow \ref{squarefree2}.$
Let $x,y \in X, x \neq y.$ By assumption $r(xx)=xx,$
so $\Lcal_x(x)=x \neq \Lcal_x(y).$
$\ref{squarefree2} \Longrightarrow \ref{squarefree1}$.
Let $x\in X$, clearly thre is an equality of sets
\[
\{\Lcal _x(y)\mid y \in X, y\neq x \}= X\setminus \{x\}
\]
so $\Lcal _x (x)= x.$
Similarly $\Rcal _x(x) = x,$ thus
$r(xx) = xx$.

For the following lemmas we assume the hypothesis of the theorem.
\begin{lemma}
\label{xy} If $(X, r)$ is nondegenerate and square-free, then
$r(xy)\neq xy$ if and only if $x \neq y$.
\end{lemma}
\begin{proof}
The statement of the lemma  follows immediately from B. and from
the equation $r(xy)=\Lcal _x(y)\Rcal _y(x)$
\end{proof}

\begin{lemma}
\label{preCCLemma} If $(X, r)$ is a non-degenerate and square-free
solution of the Yang-Baxter equation  (not necessarily
involutive), then the following conditions hold in $S$:

\begin{equation}
\label{preCC1} [yx=x^{\prime}y^{\prime}, x\neq y] \Longrightarrow
[yx^{\prime}=x^{\prime
\prime}y^{\prime},y^{\prime}x=x^{\prime}y^{\prime \prime}],
\end{equation}
for some $x^{\prime \prime}, y^{\prime \prime}  \in X.$

Furthermore, there are equalities:
\begin{equation}
\label{preCC2} yxx = x^{\prime}x^{\prime}y^{\prime
\prime},\text{and} \
 yyx=x^{\prime \prime}y^{\prime}y^{\prime}.
\end{equation}
\end{lemma}
\begin{proof}
Let $x\neq y$ and let $yx=x^{\prime}y^{\prime},$ or equivalently,
 $r(yx)=x^{\prime}y^{\prime}$. It follows from \ref{xy} that $yx
\neq x^{\prime}y^{\prime}$, as monomials in the free semigroup
$\langle X \rangle$ Assume that
\begin{equation}
\label{3}
r(yx^{\prime})=x^{\prime \prime}y^{\prime \prime}.
\end{equation}
Now consider the "Yang-Baxter diagram"
\begin{equation}
\label{diagram1}
\begin{CD}
yyx  @>r\times id_X>> yyx\\
@V  id_X\times r VV @VV id_X\times r V\\
yx^{\prime}y^{\prime}@. yx^{\prime}y^{\prime}\\
@V r\times id_X VV @VV r\times id_X V\\
x^{\prime \prime}y^{\prime \prime}y^{\prime}  @> id_X\times r >> x^{\prime \prime}y^{\prime \prime}y^{\prime}\end{CD}
\end{equation}
It follows then that
$r(y^{\prime \prime}y^{\prime}) =y^{\prime \prime}y^{\prime}$,
which, since $r$ is square-free, is possible only if
$y^{\prime\prime}=y^{\prime}.$
We have shown that
\begin{equation}
\label{eqc1}
(yx=x^{\prime}y^{\prime}) \Longrightarrow
(yx^{\prime}=x^{\prime \prime}y^{\prime}).
\end{equation}
($y^{\prime}=y$ is possible). Note
that $x^{\prime \prime}\neq y, y^{\prime}$.

Similarly, we prove that
\begin{equation}
\label{eqc2}
(yx=x^{\prime}y^{\prime}) \Longrightarrow
(y^{\prime}x=x^{\prime}y^{\prime \prime}).
\end{equation}
for some appropriate
$y^{\prime \prime}\in X$

The equality $yyx = x^{\prime \prime}y^{\prime}y^{\prime}$
in $S$ also follows from the diagram \ref{diagram1}.
\end{proof}
The valididy of  conditions C and D can be deduced from
the following lemma. Note
that in the hypothesis of the lemma we do not assume that $(X,r)$
is a solution.
\begin{lemma}
\label{LE} i) $(X, r)$ satisfies the weak cyclic condition
\ref{CCdef}.\ref{weakCC} if and only if $r$ is non-degenerate and
satisfies condition (\ref{preCC1}).

 ii) Suppose $(X, r)$ satisfies the weak cyclic
condition. Then  $r$ is involutive if and only if for every pair
$y,x\in X$  one has $\Lcal _y^x= (\Rcal _y^x)^{-1}$.
\end{lemma}
\begin{proof}

Clearly, the weak cyclic condition \ref{CCdef}.\ref{weakCC}
implies (\ref{preCC1}) and $r$ non-degenerate. Assume now that $r$
is non-degenerate and condition (\ref{preCC1}) holds.

Suppose $y,x\in X,$ $y\neq x$, and  $r(yx) =x^{\prime}y^{\prime}$
($x^{\prime} = x$, or $y^{\prime} = y$ are possible.) We denote
$x_1=x,$ $x_2 = x^{\prime},$ and apply \ref{preCC1} successively
to obtain a sequence of pairwise distinct elements $x_1, \cdots ,
x_m\in X$, such that
\begin{equation}
\label{eqc4}
r(yx_i)=x_{i+1}y^{\prime}, \ \text{for} \ 1 \leq i \leq m-1,
\text{and} \ r(yx_m)=x_1y^{\prime}.
\end{equation}
Similarly, (after an appropriate re-numeration) we obtain $y_1=y,
y_2, \cdots ,  y_k = y^{\prime} \in X,$ such that
\begin{equation}
\label{eqc5} r(y_jx_1)=x_2y_{j-1}, \  \text{for} \  2 \leq j \leq
k, \ \text{and} \   r(y_1x_1)=x_2y_m.
\end{equation}
We claim that
\begin{equation}
\label{eqc6} r(y_jx_i)=x_{i+1}y_{j-1}, \ \text{for} \ 1 \leq i
\leq m, 1\leq j \leq k,
\end{equation}
where $x_{m+1}:=x_1$, $y_0:= y_m$. We prove \ref{eqc6} by
induction on $j$.

Step 1. $j=1$. Clearly \ref{eqc4}, with $y_k = y^{\prime}$,
give the base for the induction.
Assume \ref{eqc6} is satisfied for all $j, 1 \leq j \leq j_0-1.$
We shall prove  \ref{eqc6} for $j=j_0$, $1\leq i \leq m-1$,
using induction on $i$ .
The base of the induction:
\begin{equation}
\label{eqc7}
r(y_{j_0}x_1)=x_2y_{j_0-1}.
\end{equation}
follows from \ref{eqc5}.
Assume now \ref{eqc6} is true for all $i < i_0$.
In particular,
\begin{equation}
\label{eqc8}
r(y_{j_0}x_{i_0-1})=x_{i_0}y_{j_0-1}.
\end{equation}
Then by  (\ref{preCC1}) one has:
\begin{equation}
\label{eqc9}
r(y_{j_0}x_{i_0})=t y_{j_0-1},\
\text{for some} \ t \in X.
\end{equation}
we apply  (\ref{preCC1}) again and obtain
\begin{equation}
\label{eqc10}
r(y_{j_0-1}x_{i_0})=t z,
 \end{equation}
for some $z \in X$. It follows from the inductive assumption that:
\begin{equation}
\label{eqc11}
r(y_{j_0-1}x_{i_0})=x_{i_0+1}y_{j_0-2},
\end{equation}
which together with  \ref{eqc10} gives $t=x_{i_{0}+1}$ thus
$r(y_{j_{0}}x_{i_{0}})=x_{i_{0}+1} y_{j_{0}-1}$. We have proved
that \ref{eqc6} holds for all $i$,  $1 \leq i \leq m$, and $j=
j_0,$ which verifies \ref{eqc6}. This proves i).

We set $\Lcal _y^x =(x_1, \cdots, x_m) \in Sym(X),$ and $(\Rcal
_x^y)^{-1}=(y_1, \cdots , y_k) \in Sym(X)$. Consider the
permutation
\[
\sigma _{y,x} = (x_1, \cdots ,  x_m)(y_1, \cdots , y_k).
\]
Clearly,
\begin{equation}
\label{eqc12} r(y_jx_i) = \sigma _{y,x}(x_i)\sigma
_{y,x}^{-1}(y_j).
\end{equation}

Assume now that  $r$ is involutive, and apply $r$ to \ref{eqc6} to
obtain $r(x_{i+1}y_{j-1})= y_jx_i$. This implies for $1 \leq i
\leq m$ and $1 \leq j \leq k:$
\begin{equation}
\label{eqc13} \Lcal _x^y=\Lcal _{x_i}^y=(y_1, \cdots y_k) = (\Rcal
_x^y) ^{-1} = (\Rcal _{x_i}^y)^{-1},
\end{equation}
\begin{equation}
\label{eqc14} \Lcal _y^x=\Lcal _{y_j}^x=(x_1, \cdots x_m) = (\Rcal
_y^x)^{-1} = (\Rcal _{y_j}^x)^{-1}.
\end{equation}

Conversely, \ref{eqc13} and \ref{eqc14} imply that $\sigma
_{y,x}=\sigma _{x,y}$ therefore  $r$ is involutive. This proves
the lemma, and completes the proof of the theorem.
\end{proof}
\end{proof}
\begin{remark}
\label{commentsCC} Let $(X,r)$ be an arbitrary square-free
non-degenerate solution (not necessarily involutive).
Consider the left and the right actions of $G$ on $X$, see
\ref{actions},
extending the assignment $y\rightarrow \Lcal _y$,  respectively,
$x \rightarrow \Rcal _x,$ where $\Lcal _y, \Rcal _x \in Sym(X)$
are the permutations defined via $r(yx) = \Lcal _y(x) \Rcal
_x(y).$ Since each permutation has a presentation as a product of
disjoint cycles in $Sym(X)$, (unique up to commutation of
multiples) we obtain that the cycle $\Lcal _{y}^x=(x_1, \cdots
x_m),$ $(x_1=x)$ occurs as a multiple in such a presentation of
$\Lcal_y$ and the cycle $\Rcal _x^y = (y_1, \cdots y_k)$ is a
multiple of the corresponding presentation for $\Rcal _x.$ The
surprising part is that each pair $y_j, x_i$ with $1\leq j \leq k$
and $1 \leq i \leq m,$ produces the same pair of cycles: $\Lcal
_{y_j}^{x_i}=\Lcal _{y}^x=(x_1, \cdots x_m)$, and $\Rcal
_{x_i}^{y_j}= \Rcal _{x}^{y}= (y_1, \cdots , y_k)$. Therefore
although in general $\Lcal _{y_j}\neq \Lcal _{y}$, each
permutation $\Lcal _{y_j}$ , $1 \leq j \leq k$ contains the same
cycle $(x_1,\cdots , x_m)$ in its presentation as products of
disjoint cycles in $Sym(X)$. Analogously, the cycle $(y_1, \cdots
, y_k)$ participates in the presentation of each $\Rcal _{x_i}$,
$1 \leq i \leq m,$ as a product of disjoint cycles. We do not know
how  the cycles $L_y^x$ and $R_y^x$,  are related to each other,
in the general (non-involutive) case of square-free non-degenerate
solutions, besides the obvious property, that each of them
contains $x$, see example \ref{exnoninvolutive}. In the case of
involutive solutions $(X, r)$ there is a "symmetry" $\Lcal
_y^x=(\Rcal _y^x)^{-1}= (x_1, \cdots , x_m)$ for each pair $y\neq
x, y, x \in X$.
\end{remark}

\begin{notation}
\label{leftactionnotation} To avoid complicated expressions, sometimes
we
shall use also  the notation ${}^xy=\Lcal _x(y)$ and
$y^x= \Rcal_x(y)$.
\end{notation}
The following corollary is a "translation" of the cyclic
condition in the new notation. It can be extracted
from a more general result in \cite{TShahn}.
\begin{corollary}
\label{axiomformulas}
Let $r: X\times X \r  X\times X $ be a non-degenerate
involutive bijection. Consider the following conditions:
\begin{enumerate}
\item
\label{axiomformulas1}
\[
{}^xx=x \ \text{for every}\ x\in X.
\]
\item
\label{axiomformulas2}
\[
r(x,x)=(x,x) \ \text{for every}\ x\in X.
\]
\item
\label{axiomformulas3}
$(X, r)$ satisfies the cyclic condition.
\item
\label{axiomformulas4}
For every $x, y \in X$ there are equalities:
\begin{equation}
\label{axiomformulas5}
\ {}^{({}^xy)}x={}^yx; \ ({}^yx)^y= {}^y(x^y)= x.
\end{equation}
\end{enumerate}
Then the following is true:

a)  Conditions \ref{axiomformulas1} and \ref{axiomformulas2}
are equivalent.

b)  Conditions \ref{axiomformulas3} and \ref{axiomformulas4}
are equivalent.
\end{corollary}

\begin{convention}
In the rest of the paper we shall consider only
involutive non-degenerate square-free solutions $(X, r)$ of
the Yang-Baxter equation; they will be briefly
called \emph{square-free solutions}.
\end{convention}
Let $x\in X$. Clearly, for $t \in X$ the cycle $\Lcal _x^t$ is of
length one if and only if $xt=tx^{\prime}$.
\begin{notation}
\label{Gcal}
We denote by $\Gcal _{L}=\Gcal_{L}(X, r)$ the image of $G(X, r)$
under the group homomorphism $\Lcal: G \longrightarrow Sym(X),$
which  extends the assignment $x \longrightarrow \Lcal _x$. $\Gcal
_{R}=\Gcal_{R}(X, r)$ denotes the image of $G(X, r)$ under the
group homomorphism $\Rcal: G \longrightarrow Sym(X),$ which
extends the assignment $x \longrightarrow \Rcal _x$.
\end{notation}
\begin{lemma}
\label{Lx} Let $(X, r)$ be a square-free solution, $\Lcal _x,$ and
$\Rcal _x$ be the left and right components of $r$, which are
extended to a left, respectively right action of $G( X, r)$ on $X$.
Then
\begin{enumerate}
\item The permutation $\Lcal _x$ is presented as  a product of
disjoint cycles in $Sym(X)$ via the equality:
\begin{equation}
\label{Lx1} \Lcal _x = \Lcal _x^{t_1} \Lcal _x^{t_2}\cdots \Lcal
_x^{t_s}
\end{equation}
where  $t_1, \cdots ,t_s$ are representatives of all disjoint
orbits of $\Lcal _x$ in $X$. \item The permutations $\Lcal _x$ and
$\Rcal _x$ satisfy the equality:
\begin{equation}
\label{Lx2} \Rcal _x = (\Lcal _x)^{-1}.
\end{equation}
Furthermore, the two permutation groups determined by the left and
right action of $G(X, r)$ on $X$ coincide:
\[ \Gcal _{R}=\Gcal_{L}.
\]
\item The assignment $x \longrightarrow \Lcal _x$, $x \in X$,
determines the solution $r$ uniquely, via the formula:
\[
r(x,y) =\Lcal _x(y)(\Lcal _y)^{-1}(x).
\]
\end{enumerate}
\end{lemma}


To each solution we associate an invariant integer number $M=M(X,
r)$ defined as follows.
\begin{definition}
\label{M}
\begin{enumerate}
\item For every $x\in X$ we denote by $M_x$ the order of the
permutation $\Lcal _x$ in $Sym(X)$, i.e. (in the notation of
\ref{Lx}) the least common multiple of the lengths of the cycles
$\Lcal _x^{t_i},$ $1 \leq i \leq s$. \item By $M = M(X, r)$ we
denote the least common multiple of all $M_x$, where $x\in X$, and
call $M$ \emph{the cyclic degree} of the solution $(X ,r)$.
\end{enumerate}
\end{definition}
\begin{lemma} Suppose $ax=ya^{\prime},$ for some  $ x, y, a, a^{\prime} \in X$
Then $M_x=M_y.$
\end{lemma}
\begin{proof}
It will be enough to show that the length $k$ of each cycle $\Lcal_x^\xi$ occurring
in $\Lcal_x$ divides $M_y.$
\end{proof}
\begin{proposition}
\label{Mxorbits} Assume $x, y \in X,$ and  $O_G(x)= O_G(y).$
Then $M_y = M_x.$
\end{proposition}
\begin{corollary}
Suppose $M_x \neq M_y$, for some $x, y \in X.$ Then $G$ acts
non-transitively on $X$, and $X$ is decomposable into a disjoint
union of two $r$-invariant subsets.
\end{corollary}
The following proposition  follows easily from the cyclic
condition, and \ref{M}.
\begin{proposition}
\label{newrelations} Let  $(X, r)$ be a square-free
solution of cyclic degree $M$. Let  $p$, and $q$
be arbitrary natural numbers.  Suppose $y,x \in X,$ $y\neq x$, and
let $k, m$ be the natural numbers defined in \ref{CCdef}. Let
$M_x$,denote the order of $\Lcal_x$.
Then the
following equalities hold.
\begin{equation}
\label{msmall} y^mx=xy_k^m.
\end{equation}
\begin{equation}
\label{newrelationspq} y^px^q= (x^{\prime})^{q}(y^{\prime})^{p},
\text{where} \ x^{\prime}= (\Lcal _y)^p(x),\ \text{and} \
y^{\prime}= (\Lcal _x)^{-q}(y).
\end{equation}
\begin{equation}
\label{newrelationsMx} x^{M_x} y= y (x_m)^{M_x}.
\end{equation}
\begin{equation}
\label{MM}
x^My^M=y^Mx^M.
\end{equation}
\end{proposition}
The next corollary follows immediately from \ref{MM}.
\begin{corollary}
\label{center} Let $(X, r)$ be a square-free solution, then the
center of the Yang-Baxter algebra $\Acal(k,X,r)$ contains all
symmetric functions in $x_1^M, x_2^M, \cdots, x_n^M.$
\end{corollary}
\begin{corollary}
\label{abeliansubgroup} Let $(X, r)$ be a square-free solution.
Then the group $A= gr[x_1^M, \cdots , x_n^M]$ is a free abelian
subgroup of $G (X, r)$ of index $M^n$.
\end{corollary}

\section{The lattice structure of $S(X, r)$}
\label{latticestructure} In this section we show that for a
semigroup $\Scal$ of left $I$-type, the relation $\mid_l$ of left
divisibility, defined in \ref{ldivisibilitydef}, and the left
  $I$-structure  $ v: \Ucal \longrightarrow S$ , see
  \ref{Itypedef},
  are compatible,
   and prove that
  $(S, \mid_l)$ is a distributive lattice. Analogous results
  are true for semigroups with right
  $I$-structure  $ v_1: \Ucal \longrightarrow S$.
   As a corollary we obtain
  that the Yang-Baxter semigroup $S=S(X, r)$ has a structure of
   distributive lattice, induced by its left $I$-structure
   $v$. We keep the notation from the previous sections,
   In particular,
\begin{equation}
   \Ucal=[u_1, \cdots, u_n]
   \end{equation}
is the free commutative multiplicative semigroup generated by
$u_1, \cdots , u_n$, and $\langle X \rangle$ denotes the free
semigroup generated by $X.$  The definition of an $I$- structure
is given in \ref{Itypedef}.

The following result can be extracted from \cite{TM}, Theorem 1.3.
\begin{theorem}
\label{Itype1}
Let $(X, r)$ be a square-free solution, and $S=S(X, r)$ be the associated
Yang-Baxter semigroup. Then

A. There exists a unique left $I$-structure
 $v: \Ucal \rightarrow S,$
 which is inductively defined by the following
conditions:
\begin{enumerate}
\item \label{I1} $ v_1(1)=1,$ $v(u_i)=x_i,$ for $1\leq i \leq n.$
 \item \label{I2} For every $b\in \Ucal$ and every $i, 1
\leq i \leq n,$  there exists an $x_{b,i}\in X,$ such that
$v(u_ib)=x_{b,i}v(b)$. Moreover, there is an equality of sets
\begin{equation}
\{x_{b,i}\mid 1\leq i \leq n\}=\{x_1, \cdots , x_n\}.
\end{equation}
\item \label{I3} For every $b\in \Ucal$, and $1\leq i,j \leq n,$
there is a relation in $S$:
\begin{equation}
x_{u_jb, i}x_{b, j}= x_{u_ib, j}x_{b, i}.
\end{equation}
\end{enumerate}

B. There exists a unique right $I$-structure $v_1: \Ucal \rightarrow
S,$ which satisfies the following conditions:
\begin{enumerate}
\item \label{RI1} $v_1(1)=1,$ $v_1(u_i)=x_i,$ for $1\leq i \leq
n.$
\item \label{RI2} For every $b\in \Ucal$ and every $i$, $1
\leq i \leq n$, there exists an $x_{i,b}\in X,$ such that
$v(bu_i)=v(b)x_{i,b}$. Furthermore, there is an equality of sets
\begin{equation}
\{x_{i, b}\mid 1\leq i \leq n\}=\{x_1, \cdots , x_n\}.
\end{equation}
\item \label{RI3} For every $b\in \Ucal$, and $1\leq i,j \leq n,$
there is a relation in $S$:
\begin{equation}
x_{i, bu_j}x_{j, b}= x_{j, bu_i}x_{i, b}.
\end{equation}
\end{enumerate}
\end{theorem}
\begin{remark}
Suppose $\Scal,$ is a semigroup of (left) $I$-type generated by
$x_1, \cdots, x_n,$ with a left $I$-structure $v: \Ucal
\rightarrow \Scal$. Then in general, $v$ satisfies a modified
version on condition A where condition  \ref{I1} is modified to
\begin{equation}
\label{I1general}
v(u_j)=x_{i_j}, 1 \leq j \leq n, \text{}
\end{equation}
where $i_1, \cdots , i_n$ is a permutation of
$1, \cdots, n,$
and conditions \ref{I2}, and \ref{I3} are unchanged. Moreover
\ref{I1general}
determines the bijection $v$ uniquely.
Analogous statement is true for right $I$- structures.
Without loss of generality
we can consider only the special  $I$-structures $v$ and $v_1$
defined in theorem \ref{Itype1}.
\end{remark}
\begin{notation}
\label{leftandrightistructures} Throughout  this section  $\Scal$
will denote a semigroup of $I$-type generated by $x_1, \cdots,
x_n$ with a left $I$-structure $v$ and a right $I$-structure
$v_1$. We assume that $v$ and $v_1$ satisfy conditions
\ref{Itype1} A, and B, respectively.
\end{notation}
\begin{remark}
Note  that given $a\in \Ucal$, in finitely many steps one can find
effectively the monomials $v(a)$ and $v_1(a)$. In particular, it
is easy to see that for any $i, 1 \leq i \leq n$, and any positive
integer $k$ there are equalities $v(u_i^k)= v_1(u_i^k)=x_i^k$. In
general, for  a monomial $u \in \Ucal$ there might be  inequality
$v(u) \neq v_1(u)$ (as elements of $\Scal$), see \ref{vv1}.
\end{remark}
We study first the properties of the relations  "$\mid_l$"-
\emph{divisibility with respect to left multiplication} or
shortly-
 \emph{left divisibility} and "$\mid_r$"- \emph{right divisibility}
on $\Scal$, defined as
\begin{equation}
a \mid_l b, \text{if there exists a} \
 c\in \Scal, \text{such that} \  b = ca.
\end{equation}
\begin{equation}
a \mid_r b, \text{if there exists a}  \  d\in \Scal, \text{such
that} \ b = ad.
\end{equation}
The following lemma shows that the left $I$-structure $v$ is
compatible with the left divisibility.
\begin{lemma}
\label{divisibility} $\mid_l$ is a partial order on $\Scal$,
compatible with the left multiplication. Furthermore, this order
is compatible with the left  $I$-structure $v$. More precisely,
the following two conditions hold:

a) If $a\mid b \in \Ucal$  (i.e. $b=ca$ is an equality in $\Ucal$)
then $v(a)\mid_l v(b)$;

b) Conversely, let $a, b, c \in \Scal$ satisfy $b=ca.$ Let $a_0,
b_0$ be the unique elements of $\Ucal$ which satisfy  $v(a_0)=a$
and $v(b_0)=b$. Then $b_0 = c_0a_0,$ for some $c_0\in \Ucal.$
\end{lemma}
\begin{proof}
First we show that $\mid_l$ is an ordering on $\Scal$ as a set.
Clearly, $a\mid_l a$ for every $a \in S.$ It is known that each
semigroup $\Scal$ of $I$-type is with cancellation low, see
\cite{TM}. It follows then that $a \mid_l b$ and $b \mid_l a$
imply $a=b$. The transitiveness follows at once from the
definition of $\mid_l$.

Next we prove a). Assume $b=ca$, for  $a, b, c\in \Ucal.$ We use
induction on the length $\mid c\mid$ of $c$ to find a monomial
$c^{\prime}\in \Scal,$ such that $v(b)=c^{\prime}v(a).$ If $c=
u_i,$ then by the definition of $v$ we have
$v(b)=v(u_ia)=x_{a,i}v(a).$ Assume that the statement of the
proposition is true for all $c$ of length $\leq m $. Let $b = ca,$
where $\mid c\mid = m+1.$ Then $c=u_id,$ where $1 \leq i \leq n,$
and $\mid d \mid = m.$ We have $v(b)= v(u_ida)= x_{da,i}v(da).$ By
the inductive assumption $v(da) = d^{\prime}v(a),$ so $v(b)=
x_{da,i}d^{\prime}v(a),$ which proves a). Assume now that $a,b \in
\Scal$ , and $b=ca$, for a $c\in \Scal$. By
definition,$v$ is a bijection,
so there are unique $a_0$ and $b_0$ in $\Ucal$, such
that $v(a_0)=a,$ and $v(b_0)=b.$ We have to find a $c_0\in \Ucal,$
such that $b_0=c_0a_0.$ We show this again by induction on the
length $\mid c\mid$ of $c$. If  $\mid c\mid=1$, then $c=x_i\in X$.
It follows from \ref{Itype1} that there is an equality of sets
\begin{equation}
\label{sets}
\{v(u_1a_0), \cdots , v(u_na_0)\} = \{x_1v(a_0), \cdots , x_nv(a_0)\}.
\end{equation}
Clearly, then there exists a $j$, such that
$v(u_ja_0)=x_iv(a_0)=x_ia=b$.
This gives $b_0=u_ja_0$.
Assume b) is true for all $c\in S$ with length $\mid c\mid \leq k.$
Let $b=ca$, where  $\mid c\mid = k+1$. Then $c=xd,$ for some $x\in X$
and $\mid d \mid =k$.
It follows from the inductive assumption that there is a
$d_0\in \Ucal,$
such that
\begin{equation}
\label{da}
v(d_0a_0)=dv(a)
\end{equation}
In addition an equality of sets similar to \ref{sets} shows that
there exists  an $u_j$, such that $v(u_jd_0a_0)=xv(d_0a_0).$ The
last equality together with \ref{da}  gives
$v(u_jd_0a_0)=xv(d_0a_0)=xda=ca,$ so $c_0=u_jd_0$ satisfies the
desired equality $b_0=c_0a_0.$
\end{proof}
An analogous statement is true for the right $I$-structure $v_1$.
\begin{lemma}
\label{lcm} Let $a,b \in \Scal$. a)  There exist a uniquely
determined least common multiple
of $a$ and $b$,
with respect to $\mid_l$,
that is a monomial $w$ of minimal length, such that $w=w_1a=w_2b$,
for some $w_1, w_2 \in S$. b) There exist a uniquely determined
least common multiple,
of $a$ and $b$, with respect to $\mid_r$, that is a monomial
$w^{\prime}$  of minimal length, such that
$w^{\prime}=aw_1^{\prime}=bw_2^{\prime}$, for some $w_1^{\prime},
w_2^{\prime}\in S.$
\end{lemma}
\begin{proof}
The map  $v$ is bijective, so $a=v(a_0),$ and $b=v(b_0)$, for some
uniquely determined $a_0$ and $b_0$ in $\Ucal.$ Let $w_0$ be the
least common multiple $a_0\sqcup b_0$ of $a_0$ and $b_0$ in
$\Ucal.$ It follows from \ref{divisibility}  that $v(w_0) = \xi
v(a_0)=\eta v(b_0)$. Thus $w=v(w_0)$ satisfies
\begin{equation}
\label{w}
w=\xi a= \eta b.
\end{equation}
is a common multiple of $a$ and $b$ (with respect to $\mid_l$).
That $w$ is of minimal possible length among the monomials
satisfying \ref{w} follows from \ref{divisibility}.
 This proves a). An analogous argument proves b).
\end{proof}
\begin{notation}
By $a\sqcup b$ we denote the least common multiple of $a$ and $b$
with respect to $\mid_l$. $a\vee b$ denotes the the least common
multiple of $a$ and $b$ with respect to $\mid_r.$
\end{notation}
\begin{lemma}
\label{V} Let $v,$ $v_1$ be the left and the right I-structures on
$\Scal$, defined in \ref{leftandrightistructures}. Then a) $v$ is
a lattice isomorphism for $(\Ucal, \mid)$ and $(\Scal, \mid_l)$;
b) $v_1$ is a lattice isomorphism for $(\Ucal, \mid)$ and $(\Scal,
\mid_r)$.
\end{lemma}
\begin{definition}
\label{heads} Let $u\in \Scal$. We say that $h\in X$ is a
\emph{head} of $u$ (as an element of  $\Scal$), if  $u$ can be
presented as $u=hu^{\prime}$, for some $u^{\prime} \in X.$ The
element $t\in X$ is called \emph{a tail} of $u$ (in $\Scal$) if
$u=u^{\prime\prime}t$ is an equality in $\Scal$, for some
$u^{\prime\prime}\in \Scal.$
\end{definition}
Note that a monomial may have more than one heads (respectively
tails).
\begin{example}
\label{tails} The relation $(xy=y^{\prime}x^{\prime}) \in \Re$
implies that the heads of $xy$ are $x$ and $y^{\prime}$, and
 its tails are $y$ and $x^{\prime}.$
Furthermore,  $xy=x\vee y^{\prime}= y\sqcup x^{\prime}$.
\end{example}
\begin{example}
\label{vv1} Consider the YB semigroup $S=\langle X; \Re \rangle$,
where $X=\{x_1,x_2, x_3, x_4\}$ and the set of relations is
\[
x_4x_1=x_2x_3, x_4x_2=x_1x_3, x_3x_1=x_2x_4, x_3x_2=x_1x_4,
x_1x_2=x_2x_1, x_3x_4=x_4x_3.
\]
Then
\[
v(u_2u_4)= x_1x_4=x_3x_2 = v_1(u_1u_3),
v_1(u_2u_4)=x_4x_1=x_2x_3=v(u_1u_3).
\]
\[
v(u_2^2u_4)=x_3x_2^2=x_1x_4x_2= x_1^2x_4= v_1(u_1^2u_3).
\]
\[
v_1(u_2^2u_4)=x_2^2x_4=x_2x_3x_1=x_4x_1^2.
\]
Clearly, $v(u_2^2u_4)\neq v_1(u_2^2u_4)$ as elements of $S$.
In fact,  $v(u_2^2u_4)=v_1(u_1^2u_3)$.
For $w=x_1^2x_4$ there are equalities in $S$
\[
w=x_2^2\sqcup x_4 =x_1^2 \vee x_3.
\]
\end{example}
\begin{remark}
In general, for $w\in \Ucal$ there might be an inequality
 $v(w)\neq v_1(w)$, and
 it is not true that $a\sqcup b=  a\vee b$, cf. \ref{vv1}.
Still for the special monomial
\begin{equation}
W_0=u_1u_2\cdots u_n
\end{equation}
 one has
\begin{equation}
 v(W_0)= v_1(W_0) = x_1\sqcup
x_2\sqcup \cdots \sqcup x_n = x_1\vee x_2\vee \cdots \vee x_n.
\end{equation}
\end{remark}
\begin{lemma}
Let $w_0\in \Ucal.$ Suppose
$w_0= u_{i_1}^{\alpha_1}u_{i_2}^{\alpha_2}\cdots u_{i_k}^{\alpha_k}$,
where  $1 \leq i_1<i_2< \cdots i_k \leq n,$ and all
$\alpha_j$ are positive integers.
Then

a) $v(w_0)= x_{i_1}^{\alpha_1}\sqcup x_{i_2}^{\alpha_2}\sqcup \cdots
\sqcup x_{i_k}^{\alpha_k}$;

b) $v_1(w_0)= x_{i_1}^{\alpha_1}\vee x_{i_2}^{\alpha_2}\vee \cdots
\vee x_{i_k}^{\alpha_k}$;
\end{lemma}
\begin{proposition}
\label{lattice} Let $\Scal$ be a semigroup of I-type, let  $v$ and
$v_1$ be the left and right structures on $\Scal$ as in
\ref{leftandrightistructures}. Then following conditions hold.
\begin{enumerate}
\item \label{lattice1} $(\Scal, \mid_l)$ is a distributive
lattice. More precisely,  any monomial $w\in \Scal$ has a
unique presentation as $w = x_1^{\alpha_1}\sqcup x_2^{\alpha_2}
\sqcup\cdots \sqcup x_n^{\alpha_n},$ where
$\alpha_i$ is a uniquely determined nonnegative integer
for each $i, 1 \leq i \leq n$.
In particular, for each $i,$ with $\alpha_i\geq 1,$ there is an
equality $w = w_ix_i^{\alpha_i}$, where $w_i\in \Scal,$ and $x_i$
does not occur as a tail of $w_i.$
\item The properties of the lattice $(\Scal, \mid_r)$
are analogous. In particular, every element $w\in \Scal$ has a
unique presentation as $w = x_1^{\beta_1}\vee x_2^{\beta_2}\vee
\cdots \vee x_n^{\beta_n}$, where all $\beta_i$ are nonnegative
integers. Moreover \item \label{W0} The following are equalities
in $\Scal$:
\[
W_0 = v_1(u_1u_2 \cdots u_n) =x_1\vee x_2\vee \cdots \vee x_n=
x_1\sqcup x_2\sqcup \cdots \sqcup x_n= v(u_1u_2 \cdots u_n)
\]
\end{enumerate}
\end{proposition}
\begin{proof}
It is well known that $\Ucal$ is a distributive lattice with
respect to the order of divisibility, $a\mid b.$ In particular, every element $a
\in \Ucal$ has a unique presentation $a= u_1^{k_1}u_2^{k_2}\cdots
u_n^{k_n},$ where $k_1, \cdots , k_n$ are nonnegative integers,
and $a = u_1^{k_1}\sqcup u_2^{k_2}\sqcup \cdots \sqcup u_n^{k_n}$
($v\sqcup w$ denotes the least common multiple of $v, w$ in
$\Ucal$). Lemma \ref{V} implies condition (\ref{lattice1}).
One can  show using induction on $k$  that a monomial of the
shape $u_{i_1}u_{i_2}\cdots u_{i_k},$ where all $u_{i_j}$ are
pairwise distinct, has exactly $k$ different heads and $k$
distinct tails. Therefore the monomial $W_0 =v(u_1u_2\cdots u_n)$
has exactly $n$ distinct heads (respectively, $n$ distinct tails)
so the set of heads for $W_0$ coincides with $X$.
\end{proof}
\section{Unions of  solutions and matched pairs of groups}
\label{matchedpairs}
In this section we briefly recall some definitions and properties
of unions of solutions. We also state  a recent result from \cite{TShahn},
 in which matched pairs approach is used to describe
extensions of solutions.
\begin{definition} \cite{ESS}
\label{defdisjointunion}Let $(X, r)$ be a solution.
A subset $Y\subseteq X$ is \emph{$r$-invariant}, if $r$ restricts
to a bijection $r_Y: Y\times Y \longrightarrow Y\times Y.$
$(X, r)$ is \emph{decomposable} if it can be presented as a union
of two non-empty disjoint $r$-invariant subsets.
A solution $(Z, r)$ is  \emph{ a union} of the solutions $(X, r_X)$
and $(Y, r_Y),$
if  if $X\bigcap Y= \emptyset$, $Z =
X\bigcup Y,$ as a set, and the bijection $r$ extends $r_X,$ and $r_Y$.
\end{definition}
Clearly,  $(Z, r)$ is a union of two nonempty  $r$-invariant subsets, if
and only if it is decomposable.
\begin{remark}\cite{ESS}
\label{invariant} Suppose the solution $(Z, r)$  is a
union of $(X, r_X)$ and $(Y, r_Y)$.
 Then the
map $r$ induces bijections
\[X\times Y\rightarrow Y\times X,
\text{and} \  Y\times X \rightarrow X\times Y.
\]
\end{remark}

Note that a (disjoint) union $(Z,r)$ of two square-free solutions
$(X, r_X),$ and $(Y, r_Y)$ is also a square-free solution. The
cyclic condition implies then that   for every $z \in Z$, there is
an equality $\Rcal_z = \Lcal_z^{-1}.$ Therefore the equality
$r(x,y) = (\Lcal_{x\mid Y}(y),  \Lcal_{y\mid X}^{-1}(x))$ defines
a left action of the groups $G(X,r_X)$ on the set $Y$ and a left
action of the group $G(Y, r_Y)$ on the set $ X$. Furthermore for
every $z\in Z$ there is an equality of permutations in $Sym(Z)$:
$\Lcal_z= \Lcal_{z\mid X} \Lcal_{z\mid Y}$
The following lemma is straightforward.
\begin{lemma}
\label{orbitsL} Let $(X, r)$ be a solution. Suppose
$X_1, X_2,   \cdots, X_k$ are all disjoint
orbits of the left action of $G(X, r)$ on $X$. Then $r$ induces
solutions $(X_i, r_i),$  $ 1 \leq i \leq k,$ where each $r_i$ is
the restriction of $r$ on $X_i\times X_i.$ Furthermore, $X$ is a
disjoint union of $(X, r_i)$, $1 \leq i \leq k.$
\end{lemma}
Clearly, $(X, r)$ is decomposable if and only if $G(X, r)$
acts non-transitively on $X$.
\begin{remark} W. Rump \cite{Rump} proved that every square-free solution $(X, r)$
is decomposable.
\end{remark}
Therefore to understand the structure of a solution and also
for constructing  solutions it is essential
to study extensions of solutions.
\begin{definition} \cite{ESS}
\label{ext} Suppose $(X, r_X)$ and $(Y, r_Y)$ are (disjoint)
solutions. The set of \emph{extensions of $X$ by $Y$}, denoted by
$Ext(X, Y),$ is defined as  the set of all decomposable solutions
$Z$ which are unions of $X$ and $Y$.
\end{definition}
It is shown in \cite{ESS}, that given $(X, r_X),$ and $(Y, r_Y),$
an element $Z$ of $Ext(X, Y)$ is uniquely determined by the
function : $r_{X,Y}: X\times Y \longrightarrow Y \times X.$

The fact that every square-free solution $(Z, r)$ can be presented
as a union of two disjoint solutions $(X, r_X)$, and $(Y, r_Y)$,
where the bijective map $r: Z\times Z \longrightarrow Z\times Z$
extends the maps $r_X,$ and $r_Y,$ implies that the following
theorem covers all known constructions of solutions restricted to
the square-free case.

\begin{theorem} \cite{TShahn}
\label{theoremB}
 Let $(X, r_X)$ and $(Y, r_Y)$ be disjoint
solutions, $G_X =G(X, r_X),$ $G_Y= G(Y, r_Y)$ be the groups
associated with $(X, r_X),$ and $(Y, r_Y)$, respectively. Suppose
that $Z = X\bigcup Y$, and the bijective map $r: Z\times
Z\rightarrow Z\times Z$ is an extension of the maps $r_X$ and
$r_Y.$ Then $(Z, r)$ is a solution if and only if $(G_X, G_Y)$ is
a matched pair of groups, in the sense of Majid  \cite{Shahn}.
Moreover $(Z , r)$ is square-free if and only if $(X, r_X)$ and
$(Y, r_Y)$ are square-free solutions.
\end{theorem}

\section{The equivalence of the notions square-free set-theoretic solution of the Yang-Baxter
equation, semigroup of I type, and semigroup of skew-polynomial
type} \label{equivalentnotions}

We keep all notation and conventions from the previous sections.
As usual $(X, r)$ is a square-free solution, where $X = \{x_1,
\cdots , x_n\}$ is a finite set with $n$ elements, $S=S(X, r)$,
$G= G(X, r)$,  and  $\Acal(k, X, r)$ are the associated
Yang-Baxter semigroup, group and algebra over a field $k$, defined
in \ref{associatedobjects}.  In this  section we prove Theorem
\ref{theoremA}.

For convenience of the reader, we first recall some basic
algebraic and homological properties of $S=S(X, r)$ and $\Acal(k,
X, r)$.

\begin{theorem} \cite{TM}
\label{general} Let $X$ be a finite set of $n$ elements, $(X, r)$
be a square-free solution. Let $S=S(X, r)$, $G(X, r),$ and $\Acal
= k\langle X ;\Re(r)\rangle$ be the associated
Yang-Baxter semigroup, group,  and algebra over a field $k,$
respectively. Then
the following conditions hold.
\begin{enumerate}
\item \label{general1} The semigroup $S$ is of I-type. \item
\label{general2} $S$ is a semigroup with cancellation, and $G(X,
r)$ is its group of quotients. \item \label{general3} $S$ is
Noetherian \item \label{general4} The algebra $\Acal $ is a
Noetherian domain. \item \label{general5} The Hilbert series of
$\Acal $ is $H_{\Acal}(t) = \frac{1}{(1-t)^n}$, the same as the
Hilbert series of the commutative polynomial rings in $n$
variables over $k$. \item \label{general6} $\Acal $ is Koszul.
\item \label{general7} $\Acal $ satisfies the Auslander condition.
\item \label{general8} $\Acal $ is Cohen-Macaulay. \item
\label{general9} $\Acal $ is Artin-Schelter regular ring of global
dimension $n$. \item \label{general10} The Koszul dual $\Acal ^!$
of $A$ is a Frobenius algebra. \item \label{general11} \cite{TJO}
$\Acal $ satisfies a polynomial identity. Moreover, $S$ satisfies
a semigroup identity. \item \label{general12} $\Acal$ is catenary.
\end{enumerate}
\end{theorem}
\begin{sketchproof}
For the definition of "Cohen-Macaulay" and the "Auslander
condition" see \cite{Levasseur}. Artin-Schelter regular rings are
defined in \cite{AS}. Conditions \ref{general}.\ref{general1} till
\ref{general}.\ref{general9} can be extracted from \cite{TM} (cf.
\cite{TM}, Theorems 1.3, 1.4)

Condition \ref{general}.\ref{general11} follows from a more general result in
\cite{TJO}.  It is proved  (cf. \cite{TJO}, Theorem 3.1 and Corollary
3.2)
that if a semigroup $S$ has homogeneous defining relations, and
the semigroup algebra $k[S]$ is right Noetherian and
has finite Gelfand-Kirillov dimension, then  $k[S]$
satisfies a polynomial identity, and $S$ satisfies a semigroup identity.

Condition \ref{general}.\ref{general12} follows from \cite{Sch}.
The Koszul dual algebra $A^!$ is defined in \cite{Manin}.
Condition \ref{general}.\ref{general10} follows from the fact that
a Koszul algebra $A$ of finite global dimension is Gorenstein if
and only if $A^!$ is Frobenius, cf \cite{paul}, Proposition 5.10.
\end{sketchproof}

The following theorem proofs  Conjecture \ref{mainconjecture}
\begin{theorem}
\label{proofmainconjecture} Let $(X, r)$ be a square-free
solution, where $X$ is a finite set with $n$ elements, $n \geq 2$.
Then there exists an ordering of $X= \{x_1 < x_2 <\cdots < x_n\},$
such that the Yang-Baxter semigroup $S(X, r)$ is of
skew-polynomial type (with respect to this ordering), and the
Yang-Baxter algebra $A(k,X,r)$, over an arbitrary field $k$ is a
PBW algebra with a $k$-basis the set of ordered monomials:
\[
\mathcal{N}_0=\{ x_1^{\alpha_1}x_2^{\alpha_2}\cdots
x_n^{\alpha_n}\mid \alpha_i \geq 0, 1\leq i \leq n\}.
\]
\end{theorem}
Under the hypothesis of the theorem we first prove some lemmas.
\begin{lemma}
\label{orderingXL1} There exist an ordering on $X$, $X= \{x_1 <
x_2 <\cdots < x_n\},$ such that for any pair $x,t \in X$ the
following holds.
\begin{equation}
\label{orderingX1} (tx=x^{\prime}t^{\prime})\in \Re(X,r),
\text{and}\, \ (t> x)   \Longrightarrow (x^{\prime}< t^{\prime})
\end{equation}
\end{lemma}
\begin{proof}
We use induction on $n=\mid X\mid$. Assume that the statement of
the lemma is true for all solutions $(X, r)$, with $\mid X\mid
\leq n-1.$ It follows from a theorem of Rump, \cite{Rump}, that
every square-free solution $(X, r)$, where $X$ is a finite set, is
decomposable into a disjoint union  $X = Y\bigcup Z$ of two
nonempty $r$-invariant subsets $Y,Z$. Suppose $\mid Y\mid=k$,
$\mid Z\mid=m$, $k+m=n$. Let $r_Y$ and $r_Z$ be the restrictions
on $r$ on $Y^2$ and $Z^2,$ respectively. It follows from the
inductive assumption that there exist orderings $Y= \{y_1< \cdots
< y_k \}$, and $Z=\{z_1 < \cdots < z_m\}$, which satisfy condition
\ref{orderingX1}. We set: $y_1 < \cdots < y_k <z_1<\cdots z_m$ and
verify that this is an ordering on $X$, which satisfies
\ref{orderingX1}. Assume
\begin{equation}
tx=x^{\prime}t^{\prime}\in \Re(X, r), \mbox{ and $t>x$}.
\end{equation}
We have to show that $x^{\prime}<t^{\prime}.$ Clearly if $t,x\in
Y, $ or $t,x \in Z,$ then by the inductive assumption and by the
choice of the ordering $<$, condition \ref{orderingX1} is
satisfied. Assume now $x\in Y,$ and $t\in Z.$ (Note that the case
$t \in Y, x\in Z$ is impossible since we assume $t > x$).  The
sets $Y,$ and $Z$, are $r$-invariant, therefore by \ref{invariant}
$r$ induces a map $Z\times Y \rightarrow Y\times Z.$ In particular
$tx=x^{\prime}t^{\prime}\in \Re(X, r),$ and $t\in Z, x \in Y$,
imply that $x^{\prime}\in Y, t^{\prime}\in Z.$ Hence, by the
choice of $<$, there is an inequality $x^{\prime}< t^{\prime},$
which proves \ref{orderingX1}.
\end{proof}
\begin{lemma}
\label{orderingXL2} Suppose condition \ref{orderingX1} holds. Let
$x, t \in X$, and let
 $\Lcal_t^x=(x_1, \cdots, x_k)$,
$\Lcal_x^t= (t_1, ...., t_m)$ be their associated disjoint cycles,
see \ref{CCdef}. Then $t_1 > x_1$ implies $t_j > x_i,$ for all
$i,j$, $1 \leq i \leq k$,  $1 \leq j \leq m.$
\end{lemma}
\begin{proof}
 Using induction on $i$, we first show that
\begin{equation}
\label{eqord1} t_1 > x_i, 1 \leq i \leq k.
\end{equation}
By hypothesis $t_1 > x_1,$ which gives the base for the induction.
Assume
\begin{equation}
\label{eqord2} t_1 > x_s, \mbox{for  $1 \leq s \leq i-1$}.
\end{equation}
We claim $t_1 > x_i.$ Assume the contrary,
\begin{equation}
\label{eqord3} t_1 < x_i.
\end{equation}
Note that $t_1=x_i$ is impossible, since the cycles $\Lcal_t^x$
and $\Lcal_x^t$ are disjoint. By the cyclic condition, \ref{CCdef}
one has:
\begin{equation}
\label{eqord4} t_1x_{i-1}= x_it_m,
\end{equation}
and
\begin{equation}
\label{eqord5}
t_1x_i = \left\{  \begin{array}{ll} x_{i+1}t_m &\mbox{if $i<k$}\\
                                  x_1t_m &\mbox{if $i=k$}
                  \end{array}
         \right.
\end{equation}
In the case when $i=k,$ we obtain immediately a contradiction with
\ref{orderingX1}, since
\begin{equation}
t_1x_k=x_1t_m, \mbox{and $x_1< t_1< x_k < t_m$}.
\end{equation}
Assume now $i<k.$ Then  \ref{orderingX1},  \ref{eqord5} and the
assumption  \ref{eqord3}, imply
\begin{equation}
\label{eqord6} x_{i+1}> t_m.
\end{equation}
At the same time, the equality \ref{eqord4}
and \ref{eqord2} give
\begin{equation}
\label{eqord7} t_m > x_i.
\end{equation}
We have obtained:
\begin{equation}
\label{eqord8}
x_{i+1}> t_m > x_i > t_1 >x_1 .
\end{equation}
Induction on $j$ and analogous argument show, that for $1\leq
j\leq k-i$, the following inequalities hold:
\begin{equation}
\label{eqord9}
x_{i+j}> t_m > x_i > t_1>x_1.
\end{equation}
In particular,
\begin{equation}
\label{eqord10}
x_k > t_m > x_i > t_1>x_1.
\end{equation}
Now the equality $t_1x_k=x_1t_m$ together with \ref{eqord10} give
a contradiction with \ref{orderingX1}. We have shown that
\begin{equation}
t_1 > x_i,  \  \mbox{for all $i, 1 \leq i \leq k$}.
\end{equation}
 Induction on $j$ and analogous argument show that
\begin{equation}
t_j >x_i, \ \mbox{for all $i, 1 \leq i \leq k$}.
\end{equation}
This proves the lemma.
\end{proof}
\begin{lemma}
\label{orderingXL3} Let $(X, r)$ be a square-free solution, with
an ordering $<$ on $X$ which satisfies  \ref{orderingX1}, $S
=S(X,r)$ be the associated Yang-Baxter semigroup. Then the
following two conditions are satisfied:
\begin{enumerate}
\item
\begin{equation}
\label{orderingX2} (tx=x^{\prime}t^{\prime})\in \Re(X,r),
\mbox{and ($t> x$)} \Longrightarrow (x^{\prime}< t^{\prime}),
\mbox{and ($t > x^{\prime}$)}.
\end{equation}
\item \label{orderingX3} The relations $\Re(X, r)$ form a Groebner
basis, with respect to the degree-lexicographic ordering in the
free semigroup $\langle X\rangle$,
 induced by $<$,  or equivalently the
monomials $txu,$ where $t,x,u \in X$ and $t> x > u$ do not give
rise to new relations in $S(X, r).$
\end{enumerate}
\end{lemma}
\begin{proof}
Condition \ref{orderingX2} follows immediately from Lemma
\ref{orderingXL2}. Therefore the set of defining relations $\Re
=\Re(X,r)$ for the Yang-Baxter semigroup $S(X, r)$ satisfies the
following
\begin{equation}
\label{eqor1} (x_j x_i=x_{i^{\prime}}x_{j^{\prime}}) \in  \Re \
\mbox{and ($j >i$)} \Longrightarrow (i^{\prime} < j^{\prime}),
\mbox{and ($j>i^{\prime}$)}.
\end{equation}
We have to show that $\Re$ is Groebner basis.
It follows from the theory of Groebner bases,  that each monomial
$u\in \langle X \rangle$ has a unique normal form,
denoted by $Nor(u)$, with
respect to the so called \emph{reduced Groebner basis},
$\Re_0$ which is
uniquely determined by the set $\Re$ and, $\Re\subseteq\Re_0$.
As a set $S$ can be identified
with the set of normal monomials
\begin{equation}
\label{normalmonomials}
\mathcal{N}(S) = \{ Nor(u) \mid u \in \langle X \rangle \}.
\end{equation}
Knowing the normal monomials one can uniquely restore
the set of obstructions, i.e. the set of highest monomials
in the reduced Groebner basis, $Re_0$.
To verify the equality  $\Re = \Re_0$,
therefore $\Re$ is a Groebner basis,
it will be enough to show that
the "ambiguities" $x_kx_jx_i,$ where $n \geq k
> j > i \geq 1$, do not give rise to new relations in $S(X,
r),$ or equivalently, that
each monomial of the shape $x_ix_jx_k,$ with
$1 \leq i \leq j \leq k \leq n$ is normal, with respect to
$\Re_0$.
This will follow immediately from a stronger statement:
\begin{lemma}
Each ordered  monomial $u=x_1^{\alpha_1}x_2^{\alpha_2}\cdots
x_n^{\alpha_n}$, where $ \alpha_i \geq 0, 1\leq i \leq n$
is in normal form
with respect to the reduced Groebner basis $\Re_0$
in $\langle X\rangle$.
\end{lemma}
\begin{proof}
Each relations in $\Re$ satisfies \ref{eqor1}, so
its highest monomial is $x_jx_i,$
with $j >i$, therefore
the normal form $Nor(u)$ of each $u \in \langle X \rangle$
does not contain $x_jx_i, j>i$ as a sub word.
This shows that
\begin{equation}
S  = \mathcal{N}(S) \subseteq \mathcal{N}_0,
\end{equation}
where $\mathcal{N}_0$ is the set of ordered monomials
$\mathcal{N}_0=\{ x_1^{\alpha_1}x_2^{\alpha_2}\cdots
x_n^{\alpha_n}\mid \alpha_i \geq 0, 1\leq i \leq n\}.$

The
existence of the $I$-structure $v$ on $S(X,r)$
(by definition
$v:\Ucal\longrightarrow S$ is a bijection)
implies the
equality $\mathcal{N}(S) = \mathcal{N}_0$ .
\end{proof}
We have proved \ref{orderingXL3}.
\end{proof}
\begin{prooftheorem}
The theorem follows from Lemma \ref{orderingXL3}. Note that the
Diamond Lemma \ref{B} implies that the Yang-Baxter algebra $\Acal
= \Acal (k,X,r)$ is  PBW in the sense of Priddy \cite{priddy}, and
the set of ordered monomials $\mathcal{N}_0$ projects to  a
$k$-basis of $\Acal$ (as a $k$- vector space).
\end{prooftheorem}
\begin{prooftheoremA}
The equivalence of \ref{theoremA}.\ref{theoremA1} and
\ref{theoremA}.\ref{theoremA2} follow from \cite{TM}, Theorem 1.4.
The implications $\ref{theoremA}.\ref{theoremA1}\Longrightarrow
\ref{theoremA}.\ref{theoremA3}$, and
$\ref{theoremA}.\ref{theoremA1}\Longrightarrow
\ref{theoremA}.\ref{theoremA4}$ follow from theorem
\ref{proofmainconjecture}.  Clearly, the theory of Groebner basis
implies the equivalence of conditions
\ref{theoremA}.\ref{theoremA3} and \ref{theoremA}.\ref{theoremA4}.
Theorem  1.2, \cite{TM},  proves the implication
$\ref{theoremA}.\ref{theoremA3} \Longrightarrow
\ref{theoremA}.\ref{theoremA1}$.
\end{prooftheoremA}

\section{More about $S(X, r)$ and $G(X, r)$}
\label{s(x,r)} In this section, as usual $(X, r)$ denotes a
square-free solution, where $X$ is a finite set of $n$ elements.
We show that $G=G(X, r)$ acts by conjugation on the set $X^M =
\{x_1^M, \cdots , x_n^M \}$, where $M= M(X, r)$ is the cyclic
degree of $(X, r)$ defined in \ref{M}. We compare this action with
the left action of $G(X, r)$ on the set $X$. Next we prove that
$G(X, r)$ contains a free abelian subgroup $A$ of index $M^n$, and
prove that the quotient group $\overline{G}= G/A$ can be presented
as a product of its Sylow subgroups ((cf. \ref{Sylow}). This
implies a presentation of the group $\Gcal _{L}(X, r)$ as a
product of its Sylow subgroups. As a corollary we obtain a result
of Etingof-Schedler-Solovyev, \cite{ESS}, that the group  $G(X,r)$
is solvable.
\begin{notation}
\label{(k)}
 For any positive integer $k$
we set $X^{(k)} = \{x_1^k, \cdots x_n^k \}$. By $S^k = \langle
X^{(k)} \rangle$ we denote the submonoid of $S = S(X, r)$
generated by $X^{(k)}$. If $A,B \subset S$, then as usual,  $AB$
denotes the set of all elements $u$ of the form $u=ab,$ with $a
\in A, b \in B.$
\end{notation}
\begin{proposition}
Let $k$ be a positive integer, $X^{(k)}$ and $S^k$ as in
\ref{(k)}. Then the following conditions hold.
\begin{enumerate}
\item The map $r$ induces a map $r_k:$ $X^{(k)}\times X^{(k)}$
$\longrightarrow$ $X^{(k)}\times X^{(k)}$ such that $(X^{(k)},
r_k)$ is a square-free solution. \item $S^k$ is of I-type. \item
$S^kS^j=S^jS^k$ is an equality of sets
 in $S$, for every two positive integers $k$ and $j$.
 \end{enumerate}
\end{proposition}
It follows from \ref{newrelations} that for any pair $x, y\in X$
and $M=M(X, r)$ being the cyclic degree of the solution,
there is an equality in $S:$
\[
yx^M=x_2^My, \text{where} \ \Lcal_y(x)=x_2.
\]
This implies that $G$ acts by conjugation on the set $X^{(M)}$.
The following corollary follows easily from the existence of
the $I$-structure $v$, and  \ref{newrelations}.
\begin{corollary}
\label{groupaction} Suppose $(X, r)$ is a square-free solution.,
Then
\begin{enumerate}
\item $S(X,r)$ contains the free abelian semigroup $[x_1^M, \cdots
, x_n^M]=S^M$. \item \label{} $S(X, r)$ is left and right
Noetherian. \item The group $A= gr[x_1^M, \cdots , x_n^M]$ is a
free abelian normal subgroup of $G$ of index $M^n$.  \item The group
$G=G(X,r)$ acts by conjugation on the set $X^{(M)}$.  Moreover the
action of $A$ on $X^{(M)}$ is trivial, thus the quotient group
$\overline{G}= G/A$ acts on $X^{(M)}$ by conjugation. Clearly,
$\overline{G}$ is a finite group of order $M^n$.
\item
\label{overlineLcal}
The group $A$ is contained in the kernel $ker\Lcal$
of the homomorphism $\Lcal:G\longrightarrow Sym(X).$
Therefore there exists an epimorphism
$\overline{\Lcal}\overline{G}\longrightarrow \Gcal _{L}$,
induced by $\Lcal,$  satisfying the equality:
$\Lcal= \overline{\Lcal}\circ \nu,$ where $\nu$
is the natural epimorphism
$\nu: G \longrightarrow \overline{G}.$
\item
The order of  $\Gcal _{L}$ divides $M^n.$
\end{enumerate}
\end{corollary}
\begin{notation}
\label{equivalence1} For every $y \in X$ we denote by $O(y^M)$ the
orbit of $y^M$ under the action if $G$ on $X^{(M)}$. For $x,y \in
X$ we define an equivalence  on $X$ by setting $x \approx y$ iff
$O(x^M)=O(y^M).$ By $X(y)$ we denote the equivalence class of $y$,
$y\in X$.
\end{notation}
The lemma below follows  straightforward from the definition of
the actions of $G$ on the sets $X$ and $X^M$, and from Proposition
\ref{newrelations}.
\begin{lemma}
\label{orbits}
The following conditions hold.
\begin{enumerate}
\item There exists a one-to-one correspondence between the
G-orbits of $X^M$ and the G-orbits of $X$. More precisely for
every $\xi\in X$, there are equalities $O_G(\xi)= X(\xi)=\{ x \in
X \mid x^M \in O(\xi^M)\}$. Furthermore, the orbits $O_G(\xi)$ can
be obtained simply by acting with the "semigroup" elements of $G$:
e.g. $y \in O_G(x)$ if and only if, there exist monomials $a, b
\in S$, $a = a_1\cdots a_k $, and $b = b_1 \cdots b_k$  ( $a_i,
b_i \in X$) and elements $y_1, \cdots , y_k \in X,$ such that
there are equalities:
 \begin{equation}
\label{eqorbit}
a_kx=y_kb_k, \ a_{k-1}y_k=y_{k-1}b_{k-1},\  \cdots , \ a_1y_1=yb_1.
\end{equation}
\item If $x\in X(a)$, and $y \in X(b)$, for some $a, b \in X$ (not
necessarily $a \neq b$) then there is an equality $xy =
y^{\prime}x^{\prime},$ with $y^{\prime} \in X(b)$, $x^{\prime}\in
X(a)$. \item Each orbit $O_G(\xi), \xi\in
X$ is $r$- invariant.
\item $X$ is $r$-decomposable if and only if $\overline{G}$
does not act transitively on $X^M$. More precisely, if
$O_{\overline{G}}(\xi_i^{M}), 1 \leq i \leq k$ are all
disjoint orbits of this action  then  $X$ splits into a disjoint
union of $k$ nonempty $r$-invariant subsets: $X_1=O_G(\xi_1), \cdots ,
X_k=O_G(\xi_k)$.
\end{enumerate}
\end{lemma}






\begin{remark}
\label{orderorbit}
It is a routine fact, that the order of each orbit $O(x^M)$, $x \in X$ is a divisor of the order
$M^n$ of $\overline{G}$, see for example \cite{A1}, 6.1.
\end{remark}
A sufficient condition for $r$-decomposability of  $X$ follows
immediately from \ref{orderorbit}. As a corollary we obtain a
result  from \cite{ESS}, that every solution $(X, r)$, where $X$
is of prime order $p$ is decomposable.
\begin{corollary}
\label{p} If $M$ is not divisible by some prime divisor $p$ of $n,$
then the action of $\overline{G}$ (and of $G$) on $X^{(M)}$ is not
transitive and $X$ is a disjoint union of $k$  $r$-invariant
subsets, where $k \geq 2$ is the number of orbits in $X^{(M)}$.
\end{corollary}
\begin{corollary} \cite{ESS}
If $n = p$ is a prime number, then $X$ is a disjoint union of  two
nonemty $r$-invariant
subsets.
\end{corollary}

Next we study the relations between the cyclic degree $M=M(X, r),$
 the Sylow subgroups of $\overline{G}$
and the cyclic properties of the semigroup $S(X,r)$.
Note that
\begin{notation}
\label{q_i}
Let $M=M(X, r)$ be the cyclic degree of the solution $(X, r)$
defined in \ref{M}.
Suppose $M=p_1^{\alpha_1}p_2^{\alpha_2} \cdots p_k^{\alpha_k}$
where $p_1, \cdots , p_k$ are distinct prime numbers, and
$\alpha_1 \cdots \alpha_k $ are  positive integers.
For $i = 1, \cdots , k$, we set
\[
q_i = M / p_i^{\alpha_i},
\]
\[
S^{q_i} = \langle x_1^{q_i}, \cdots ,  x_n^{q_i}\rangle,
\]
the sub-monoid of $S$ generated by $x_1^{q_i}, \cdots ,
x_n^{q_i}$,
$1 \leq i \leq k.$
We denote  by
$\overline{S^{q_i}}$
the natural image of $S^{q_i}$ in the quotient group $\overline{G}$, and
by $\Lcal({S^{q_i}})$ the image of
$S^{q_i}$ under the homomorphism
$\Lcal: G \longrightarrow \Gcal _L \subset Sym(X),$
 defined by the left action of $G$ on $X$.
\end{notation}

Clearly, the integers $q_1, \cdots , q_k$ are pairwise coprime,
and $\overline{S^{q_i}}$ are submonoids of $\overline{G}$.

The next theorem gives a presentation of  $\overline{G}$ as a product of its Sylow subgroups.
Surprisingly it also allows to consider each element
of $\overline{G}$ as an element of the monoid $\overline{S}$.

\begin{theorem}
\label{Sylow}
The following conditions hold.
\begin{enumerate}
\item
\label{Sylow1}
For every $i$,  $1 \leq i \leq r,$ the submonoid
$\overline{S^{q_i}}$ is a subgroup of order $p_i^{n\alpha_i}$ in
$\overline{G}$.
In particular, it is a  Sylow $p_i$-subgroup
of $\overline{G}$.
\item
\label{Sylow2}
For every pair $q_i, q_j$ , $1 \leq i,j \leq r,$
there is an equality
$\overline{S^{q_i}}.\overline{S^{q_j}} =\overline{S^{q_j}}.\overline{S^{q_i}}$.
\item
\label{Sylow3}
The group $\overline{G}$ is a product of its Sylow subgroups:
$\overline{G}=\overline{S^{q_1}} \cdots \overline{S^{q_k}}.$
In particular ,  $\overline{G}=\overline{S}$.
\item
\label{Sylow4}
For each $i, 1 \leq i \leq k,$ such that $\Lcal(S^{q_i})\neq id_X$,
 $\Lcal(S^{q_i})$ is a $p_i$-Sylow subgroup of $\Gcal _L.$
\item
\label{Sylow5}
 Let $1\leq i_1, \cdots, i_s\leq k$ be all indices, for which
$\Lcal(S^{q_{i_j}})\neq \{id_X \}, 1 \leq j \leq s.$
Then the  group \ $\Gcal _L= \Gcal _L(X, r)$ is a product of
its Sylow subgroups:
\[
\Gcal _L=\Lcal({S^{q_{i_1}}}) \cdots \Lcal(S^{q_{i_s}}) .
\]
In particular, $\Gcal _L= \Lcal(S).$
\item
\label{Sylow6}
The groups $\Gcal _L$,  $\overline{G},$ and $G$ are solvable.
\end{enumerate}
\end{theorem}
\begin{proof}

Consider  $\overline{S^{q_i}}$, where $1 \leq i \leq k$. Note
first that as a finite submonoid of the group $\overline{G}$,
$\overline{S^{q_i}}$ is a subgroup of   $\overline{G}$. We claim
that the order of  $\overline{S^{q_i}}$ is exactly
$p_i^{n\alpha_i}$. The equalities
\ref{newrelationspq} imply that every element
$w$ of  $\overline{S^{q_i}}$  can be presented as
 \begin{equation}
\label{overline}
w=\overline{v((u_1^{q_i})^{\beta_1}\cdots (u_n^{q_i})^{\beta_n}}, \
\text{where}
\ 0\leq \beta_s \leq p_i^{\alpha_i}
\  \text{for all } \  s, \ 1\leq s \leq n.
\end{equation}
We set $\beta = (\beta_1, \beta_2, \cdots ,   \beta_n)$, and
$w=w(\beta),$
for the monomial $w$ determined by \ref{overline}.
 It follows from the properties of the I-structure $v$ on $S$
 and from \ref{newrelationspq}
 that each inequality  $\beta^{\prime}\neq \beta^{\prime\prime}$
implies an inequality in $\overline{S}$:
\begin{equation}
\label{im1}
w(\beta^{\prime})\neq w(\beta^{\prime\prime}).
\end{equation}
 This implies that
$\overline{S^{q_i}}$ is a group of order $(p_i^{\alpha_i})^n$ thus
a Sylow $p_i$ subgroup of $\overline{G},$ which proves \ref{Sylow1}.

Next we recall that for every pair of integers $i, j$,
$1 \leq i,j \leq k,$
and for every pair $x,y \in X$ there exist $z, t \in X,$
such that the equality
\begin{equation}
\label{im2}
x^{q_i}y^{q_j}=z^{q_j}t^{q_i}
\end{equation}
holds in $S$. This implies that
$\overline{S^{q_i}}\overline{S^{q_j}}=\overline{S^{q_j}}\overline{S^{q_i}}$
for all $i, j$, which verifies \ref{Sylow2}. Let  $S^{\prime}=\langle
S^{q_1}, \cdots S^{q_k}\rangle$ be the submonoid of $S$, generated
by $S^{q_1}, \cdots S^{q_k}$. It follows from \ref{im2} that there
is an equality
\begin{equation}
\label{im3} S^{\prime}= S^{q_1} \cdots S^{q_k}.
\end{equation}
Hence
\begin{equation}
\label{im4}
\overline{S^{\prime}}=\overline{ S^{q_1}} \cdots \overline{S^{q_k}}.
\end{equation}
is a presentation of $\overline{S^{\prime}}$ as a product of
subgroups with pairwise co-prime orders: $p_1^{n\alpha_1}, \cdots
, p_k^{n\alpha_k},$ respectively. It follows then that the order
of $ \overline{S^{\prime}}$ is exactly $p_1^{n\alpha_1}\cdots
p_k^{n\alpha_k} =M^n$, thus $\overline{G}=\overline{ S^{q_1}}
\cdots \overline{S^{q_k}}.$ This proves \ref{Sylow3}. The proof
of \ref{Sylow4}, and \ref{Sylow5} is routine.
\end{proof}
Note that, in general the Sylow subgoups $\overline{ S^{q_i}}$
might not be normal subgroups of  $\overline{G},$
as shows the following example.
\begin{example}
Let $S=\langle X; \Re \rangle$, where $X=\{x_i\mid 1\leq  i \leq
6\}\bigcup  \{y_j \mid 1 \leq j \leq 4\}$ and the relations $\Re$
are defined by the permutation
\begin{equation}
\sigma= (x_1 x_2 x_3 x_4 x_5 x_6)(y_1 y_2 y_3 y_4);
\end{equation}
as follows:
\begin{equation}
y_jx_i=\sigma(x_i)\sigma^{-1}(y_j),  \and \ x_iy_j = \sigma_1(y_j)
\sigma^{-1}(x_i) \ \text{for} \ 1\leq i \leq 6;  \  1\leq  j \leq
4;
\end{equation}
\begin{equation}
x_ix_k=\sigma^3(x_k) \sigma^{-3}(x_i), \ \text{for all}\ i\neq k
\text{(mod $3$)},\  1\leq i , k \leq 6 ;
\end{equation}
\begin{equation}
x_ix_k=x_kx_i, \ \text{for all} \  i = k  \text{(mod $3$)},\  1
\leq i,k, \leq 6
\end{equation}
\begin{equation}
y_jy_k=\sigma^2(y_k) \sigma^{-2}(y_j), \  \text{for all}\ j \neq k
\ \text{(mod $2$)}, \ 1\leq j, k \leq 4 .
\end{equation}
\begin{equation}
y_jy_k=y_ky_j, \  \text{for all}\ j = k \ \text{(mod $2$)}, \
1\leq j, k \leq 4 .
\end{equation}
It is easy to verify that the  set of relations  $\Re$ defines
naturally a square-free solution, $r$, thus $S$ is an YB
semigroup. The set of all lengths of cycles is $6,  4,  2$, thus
$M= 12=2^2.3,$ and (in the notation \ref{q_i}), $q_1=3,$ and
$q_2=4.$ Thus,  by Theorem \ref{Sylow},
$\overline{G}=\overline{S^3}\overline{S^4}$. Note that none of the
subgroups $\overline{S^3},$  $\overline{S^4}$ is normal in
$\overline{G}$
\end{example}
One can use Theorem \ref{Sylow}  to give a
straightforward proof of the $r$- decomposability of $(X, r)$  in
all cases when the cycles are not enough "dense" on $X$. More
precisely, the following corollary is true.
\begin{corollary}
\label{extreme} Suppose that there exists a prime divisor $p$ of
$n$, and an $x \in X$, such that $x$ does not belong to a cycle of
length divisible by $p.$ Then the action of $\overline{G}$ on $X$
is non-transitive, therefore
 $(X, r)$ is decomposable.
\end{corollary}
\section{Multipermutation solutions and generalized twisted unions}
\label{multipermutation}

We give a description of the generalized twisted unions of
solutions $Z = X\bigcup Y$, showing that the group $G_Y=G(Y, r_Y)$
acts as automorphisms on $X$, and all the elements $\xi$ of an
orbit $O(x)= O_{G_Y}(x)$ have the same action on $Y$ see
\ref{twistedunionsP}.  Lemma \ref{equivalenceL3} generalizes the
cyclic condition. We give a conjecture that every multipermutation
solution of level $m$ is a generalized twisted union of
multipermutation solutions of level $\leq m-1$.
 We keep the notation and conventions from the previous
sections. In particular, to avoid complicated expressions
sometimes we shall use both
notation ${}^xy=\Lcal _x(y)$ and
$y^x= \Rcal_x(y).$

\begin{definition}
\label{twistedunionsD} \cite{ESS} Let $(Z, r)$, be a disjoint
union of the solutions $(X, r_X)$, and $(Y, r_Y)$.
\begin{enumerate}
\item $(Z, r)$ is called  \emph{a twisted union} of $X$ and $Y$ if
the maps $r_{XY}: X\times Y \rightarrow Y\times X$ and $r_{YX}:
Y\times X \rightarrow X\times Y$ are defined as
\begin{equation}
 r_{XY}(x, y) = (g(y), f^{-1}(x))
 \end{equation}
  and
  \begin{equation}
  r_{YX}(y, x) = (f(x), g^{-1}(y)),
\end{equation}
where $f \in Sym(X), $ and $g \in Sym(Y)$ are fixed. \item $(Z,
r)$ is \emph{a generalized twisted union} of $X$ and $Y$ if the
map $r$ is determined by the formula:
\begin{equation}
r_{XY}(x,y)=(\Lcal_ {x\mid Y}(y), \Rcal _{y\mid X}(x)),
\end{equation}
where the permutations $\Lcal_ {x\mid Y} \in Sym(Y)$,  and  $\Rcal
_{y\mid X} \in Sym(X)$ satisfy the following condition:

(*) \emph{For every $y \in Y$ the permutation $\Lcal_ {x^y\mid Y}:
Y\rightarrow Y$ is independent of $y$, and for every $x\in X$, the
permutation $\Rcal _{{}^xy\mid X}: X \rightarrow X$ is independent
of $x$}.
\end{enumerate}
\end{definition}
\begin{notation}
When the element $\xi\in Z$ is specified we shall simply write, as
usual, $\Lcal_x(\xi),$ or ${}^x\xi$ instead of $\Lcal_{x\mid Y}(\xi),$
respectively $\Lcal_{y}(\xi),{}^y\xi$ instead of $\Lcal_{y\mid X}(\xi).$
\end{notation}
\begin{proposition}
\label{twistedunionsP} Let $(Z, r)$ be union of the disjoint solutions
 $(X, r_X),  and (Y, r_Y).$
Then $(Z, r)$ is a
generalized twisted union of $X$ and $Y$ if and only if for every
pair $x,y$, $x \in X$, $y\in Y$ the following equalities hold:
\begin{equation}
\label{twistedunionsP1} \Lcal_ {x^y\mid Y}=\Lcal _{x\mid Y} =
\Lcal_ {{}^yx\mid Y};
\end{equation}
\begin{equation}
\label{twistedunionsP2} \Lcal _{{}^xy\mid X}=\Lcal _{y\mid X} =
\Lcal_ {y^x\mid X};
\end{equation}
\end{proposition}
\begin{proof}
Note first that the  equalities \ref{twistedunionsP1} and
\ref{twistedunionsP2} imply that $(Z, r)$ is a generalized twisted
union of $X$ and $Y$.

Assume now that $(Z, r)$ is a generalized twisted union of $X$ and
$Y$.
 Let $x \in X$, $y\in Y$. We have to show
that for every $z\in Y$ there is an equality
\begin{equation}
\Lcal _{x^y}(z) =\Lcal _x (z).
\end{equation}
By  definition \ref{twistedunionsD} the map $\Lcal _{x^y\mid
Y}:Y\rightarrow Y$ is independent of $y \in Y$. Hence for every
pair $y, z \in Y$ there is an equality
\begin{equation}
\label{tu1} \Lcal _{x^y}(z) = \Lcal _{x^z}(z)
\end{equation}
By the cyclic condition in $(Z, r)$, see \ref{axiomformulas5},
one has:
\begin{equation}
\label{tu2} \Lcal _{x^z}(z)= \Lcal _x(z).
\end{equation}
Now the equations \ref{tu1} y \ref{tu2} imply
\begin{equation}
\Lcal _{x^y}(z)= \Lcal _x(z)
\end{equation}
for every $z \in Y.$
We have shown that
\begin{equation}
\label{tu3} \Lcal _{x^y\mid Y}= \Lcal _{x\mid Y}
\end{equation}
 for arbitrary $x\in X$ and $y\in Y$.
We apply this to the pair ${}^yx \in X $ and $y \in Y$ and obtain
\begin{equation}
\label{tu4} \Lcal _{({}^yx)^y\mid Y} = \Lcal _{{}^yx \mid Y}.
\end{equation}
By \ref{axiomformulas5} there is an equality:
\begin{equation}
\label{tu5} ({}^yx)^y = x,
\end{equation}
which together with \ref{tu4}, \ref{tu3} implies
 $\Lcal _{{}^yx \mid Y}= \Lcal _{x\mid Y}=\Lcal _{x^y\mid Y}.$

 This completes the proof of \ref{twistedunionsP1}.
Analogous argument proves   \ref{twistedunionsP2}.
\end{proof}
\begin{theorem}
Let $(Z, r)$ be a generalized twisted union of the solutions
$(X, r_X)$ and $(Y, r_Y)$,  and let $G_X=G(X, r_X)$, $G_Y=G(Y, r_Y)$
be the associated Yang-Baxter groups. Suppose $O_{G_Y}(\xi _1), \cdots ,
 O_{G_Y}(\xi _p)$ are the  (distinct) orbits of the
action of the group $G_Y$ on $X$, and
$O_{G_X}(\eta_1), \cdots ,  O_{G_X}(\eta_q)$ are the (distinct)
orbits of the action of $G_X$ on $Y$. Then the following
conditions are satisfied.
\begin{enumerate}
\item The assignment
\[
 x\rightarrow \Lcal _{x\mid Y}, \mbox{for all} \ x \in X
\]
extends to a group homomorphism
\[
L _X: G(X, r_X) \rightarrow Aut(Y, r_Y)
\]
\item
Let $H_X$ denotes the kernel $KerL _X$. Then
each orbit
$O_{G_Y}(\xi _i)$, is contained in the
left coset $\xi _iH_X$, i.e.
$O_{G_Y}(\xi _i)\subseteq \xi _iH_X.$
In particular, for every $x \in  O_{G_Y}(\xi _i)$ , $1\leq i \leq p$
there is an equality
\begin{equation}
\Lcal _{x\mid Y}=\Lcal _{\xi _i\mid Y}.
\end{equation}
\item
The assignment
\[
y\rightarrow \Lcal _{y\mid X}, \text{for all} \ y \in Y
\]
extends to a group homomorphism
\[
\L_Y: G(Y,r_Y) \rightarrow Aut(X,r_X).
\]
\item
Let $H_Y$ denotes the kernel $KerL _Y$. Then
$O_{G_X}(\eta _j)\subseteq \eta _jH_Y,$ for
$1\leq j\leq q.$

In particular, for every $y \in  O_{G_X}(\eta _j)$ ,  there
is an equality:
\begin{equation}
\Lcal _{y\mid X}=\Lcal _{\eta _j\mid X}.
\end{equation}
\end{enumerate}
\end{theorem}
\begin{definition}
\label{retract} \cite{ESS} Let $(X, r)$ be a square-free solution.
Define an equivalence relation on $X$ as $x \sim y $ iff
$\Lcal_x=\Lcal _y$.

Clearly , since $\Rcal _x = \Lcal _x^{-1}$, one has also $x\sim y
$ iff $\Rcal_x=\Rcal _y.$ Let $X^{\sim}= X/\sim.$ It is known, see
\cite{ESS}, that the solution  $r: X\times X \rightarrow X\times
X$ induces a bijection $r^{\sim}: X^{\sim}\times X^{\sim}
\rightarrow X^{\sim}\times X^{\sim},$ so that $( X^{\sim},
r^{\sim})$ is a solution. It is not difficult to see that this
solution is also square-free. The solution $X^{\sim}, r^{\sim})$
is called the \emph{retraction of} $(X, r)$ and is denoted by
$Ret(X, r).$ The solution is \emph{retractible} if $\sim$ is a
nontrivial equivalence relation. (or equivalently $Ret(X,r)\neq
(X,r).$ In the case when $\sim$ is the trivial equivalence on $X$,
the solution $(X,r)$ is called \emph{irretractible.}
\end{definition}
\begin{lemma}
\label{equivalenceL1} For any $x, y \in X$ the equivalence $x \sim
y$ implies $xy=yx.$
\end{lemma}
\begin{definition}
\label{retractk} Inductively, for $1 < k $ we define the
retractions of higher level as $Ret^{k}(X,
r)=Ret(Ret^{k-1}(X,r)).$

 We denote by $x^{(k)}$ the image of $x$ in
$Ret^{k}(X, r)$.  The set
\begin{equation} [x^{(k)}]:=\{ \xi \in
X\mid \xi^{(k)}=x^{(k)}\}
\end{equation}
is called \emph{the $k^{th}$ retract orbit} of $x$.
\end{definition}
\begin{definition}
\label{multipermutationalDef} \cite{ESS},  A solution $(X, r)$ is
called \emph{multipermutation solution of level} $m$ if $m$ is the
minimal nonnegative integer, such that $Ret^m(X, r)$ is finite of
order $1$.
\end{definition}
\begin{lemma}
\label{equivalenceL2} For any positive integer  $k$, and any $x\in
X$ the $k^th$ retract orbit $[x^{(k)}]$ is $r$-invariant.
Furthermore, if we denote by $\textbf{r}_{x,k}$ the corresponding
solution induced by $r$, then $([x^{(k)}], \textbf{r}_{x,k})$ is a
multipermutation
 solution of level $k$.
\end{lemma}
\begin{lemma}
\label{equivalenceL3} Let $(X, r)$ be a square-free solution. Then
the following conditions hold:
\begin{enumerate}
\item \label{equivalenceL31}  For every  $x, y, t \in X$, and  $k$
a positive integer,
\begin{equation}
y^{(k)}=t^{(k)} \Longrightarrow ({}^yx)^{(k-1)}= ({}^tx)^{(k-1)}.
\end{equation}
\item \label{equivalenceL32} For every  $x, y, t \in X$
\begin{equation}
\label{equivalenceL33}
y^{(2)}=t^{(2)} \Longrightarrow {}^yx \sim {}^tx, \ \text{in
particular}, \ {}^yt \sim t, \text{and} \ {}^ty \sim y.
\end{equation}
\end{enumerate}
\end{lemma}
\begin{proof}
We first prove   \ref{equivalenceL31}. By hypothesis,
$y^{(k)}=t^{(k)},$ or equivalently
\begin{equation}
\label{eqsym1}
 y^{(k-1)} \sim t^{(k-1)}.
\end{equation}
 Let $x \in X$. Clearly,
\begin{equation}
\label{eqsym2} yx= \xi y_1, tx=\xi _1 t_1, \  \text{for some} \
\xi, \xi _1, y_1, t_1 \in X.
\end{equation}
This implies the following equalities in $Ret^{k-1}(X, r)$
\begin{equation}
\label{eqsym3} y^{(k-1)}x^{(k-1)}= \xi^{(k-1)} y_1^{(k-1)}, \
\text{and} \
t^{(k-1)}x^{(k-1)}=\xi _1^{(k-1)} t_1^{(k-1)}.
\end{equation}
It follows then from \ref{eqsym1} that
\begin{equation}
\label{eqsym4} \xi^{(k-1)} = \xi_1^{(k-1)}, \text{or
equivalently}, \
 \xi^{(k-2)} \sim \xi_1^{(k-2)}.
\end{equation}
By \ref{eqsym2},   one has $\xi = {}^yx$, and $\xi_1 = {}^tx,$
which proves \ref{equivalenceL31}. Condition \ref{equivalenceL32}
follows straightforward from \ref{equivalenceL31}, with $k=2$, and
the cyclic condition.
\end{proof}
\begin{corollary}
\label{equivalenceC}
Let $(X, r)$ be a multipermutation solution of level $m$,
$G_X= G(X,r)$ be the associated Yang-Baxter group.
Then for every $y \in X$ one has:
\[
O_{G_X}(y)\subseteq [y^{(m-1)}],
\]
where
$O_{G_X}(y)$ is the $G_X$ orbit of $y$ in $X,$ and
 $[y^{(m-1)}]$ is the $(m-1)$-th retract orbit of $y$.
 In particular, $O_{G_X}(y)$ is a multipermutation
 solution of level at most $m-1$.
\end{corollary}
The cyclic condition, ${}^{({}^yx)}(y) = {}^xy$ is "extended" to
the class $[y^{(2)}]$ by the following lemma.
\begin{lemma}
\label{equivalenceL4}  Let $(X, r)$ be a solution. Then the
following conditions hold.
\begin{enumerate}
\item \label{equivalenceL41} For every $x \in X$, and  $z\in
[y^{(2)}]$ there is an equality
\begin{equation}
\label{equivalenceL42}
 {}^{({}^yx)}(z) = {}^xz.
\end{equation}
and
\begin{equation}
\label{equivalenceL43}
 \Lcal_{{}^yx\mid [y^{(2)}]} = \Lcal_{x\mid
[y^{(2)}]}.
\end{equation}
\item \label{equivalenceL44}
 Suppose that  $[x^{(2)}] \neq [y^{(2)}]$, and
the set $[y^{(2)}]$ is invariant under the left action of $G
([x^{(2)}], \textbf{r}_{x,2} )$, respectively, $[x^{(2)}]$ is
invariant under the left action of $G ([y^{(2)}], \textbf{r}_{y,2}
).$ Then the disjoint union $Z=[x^{(2)}]\bigcup [y^{(2)}]$ is a
generalized twisted union of $[x^{(2)}]$ and $[y^{(2)}]$.
Moreover, $(Z, r_Z)$ is a multipermutation solution of level $3$,
where $r_Z$ is the restriction of $r$ on $Z\times Z.$
\end{enumerate}
\end{lemma}
\begin{proof}
Let $x\in X$, and let $z \in y^{(2)}.$ We will show
that \ref{equivalenceL42} holds.
It follows from \ref{equivalenceL33} that
\begin{equation}
{}^yx \sim {}^zx.
\end{equation}
So, by the definition of $\sim$, and by the cyclic condition,
\begin{equation}
{}^{({}^yx)}(z) = {}^{({}^zx)}(z)= {}^xz
\end{equation}
We have shown  \ref{equivalenceL42}. Clearly, \ref{equivalenceL42}
implies \ref{equivalenceL43}. Condition
\ref{equivalenceL44} follows  easily from \ref{equivalenceL31}.
\end{proof}
\begin{corollary}
Let $(X, r)$ be a multipermutation solution of level $3$. Then
 $(X, r)$ is a generalized twisted union of multipermutation
solutions of level $\leq 2.$
\end{corollary}
\begin{example}
Let $X= \{x, x_1, \xi, \xi_1, t, t_1, \eta, \eta_1, y, y_1 \}$ and
let $r$ be determined via
\begin{equation}
\Lcal_x = \Lcal_ {x_1}=(t  t_1) (\eta \eta_1)(y y_1);
\end{equation}
\begin{equation} \Lcal_{\xi}= \Lcal_{\xi_1} = (t \eta)(t_1 \eta_1) (y
y_1);
\end{equation}
\begin{equation}
\Lcal_t=\Lcal_{t_1} = \Lcal_{\eta}= \Lcal_{\eta_1}= id_X
\end{equation}
\begin{equation}
\Lcal_y = \Lcal_{y_1}= (x \xi)(x_1 \xi_1).
\end{equation}
Then $Ret(X, r) = (X^{\sim}, r^{\sim})$, where $X^{\sim}=\{
x^{\sim}, \xi^{\sim}, t^{\sim}, y^{\sim}\}$, and  $r^{\sim}$ is
determined by $\Lcal y^{\sim}=(x^{\sim}\xi^{\sim})$,
$\Lcal_{x^{\sim}}=\Lcal_ {\xi^{\sim}} = \Lcal_
{t^{\sim}}=id_{X^{\sim}}.$ Clearly, $X^{(2)}= \{ y^{(2)},
x^{(2)}\}$, and  $Ret^{2}(X, r)$ is the trivial solution,
therefore $Ret^3(X, r)= 1.$ In this case $(X, r)$ is a
multipermutation solution of level $3$.
\end{example}

\section{Binomial solutions of the classical Yang-Baxter equation}
\label{binomialYBE}
In this section we study a particular class of solutions of the
classical Yang-Baxter equation, called \emph{binomial solutions}.
We show that there is a close relation between a class of
Artin-Schelter regular rings, which we call skew-polynomial rings
with binomial relations and the square-free binomial solutions of
the classical Yang-Baxter equation.
\begin{definition}
Let $V$ be a finite dimensional vector space over a field $k$ with
a $k$- basis $X=\{ x_1, \cdots, x_n \}$. Suppose the linear
automorphism $R:V\otimes V\longrightarrow V\otimes V$ is a
solution of the Yang-Baxter equation. We say that $R$ is \emph{a
binomial solution of the (classical) Yang-Baxter equation} or
shortly \emph{binomial solution} if the following conditions hold
\begin{enumerate}
\item for every pair $i\neq j, 1 \leq i,j \leq n,$
\begin{equation}
\label{RB} R(x_j\otimes x_i) = c_{ij} x_{i^{\prime}}\otimes
x_{j^{\prime}},
  R(x_{i^{\prime}}\otimes
x_{j^{\prime}}) = \frac{1}{c_{ij}} x_jx_i,\
\text{where}\  c_{ij} \in
k, c_{ij}\neq 0.
\end{equation}
\item $R$ is \emph{non-degenerate}, that is \emph{the associated
set-theoretic solution } $(X, r(R))$, where $r=r(R): X\times X
\longrightarrow X\times X$ is defined as
\begin{equation}
\label{associatedr}
r(x_j, x_i) = (x_{i^{\prime}}, x_{j^{\prime}}) \ \text{if} \
R(x_j\otimes x_i)
= c_{ij} x_{i^{\prime}}\otimes x_{j^{\prime}},
\end{equation}
is non-degenerate.
\end{enumerate}
We call the binomial solution $R$ \emph{square-free} if
$R(x_i\otimes x_i) = x_i\otimes x_i$, or equivalently, $(X, r)$ is
square-free.
\end{definition}
\begin{notation} By $(k, X, R)$ we shall denote a square-free
binomial solution of the classical Yang-Baxter equation.
\end{notation}
Each square-free binomial solution $(k, X, R)$ defines a quadratic
algebra $\Acal_R= \Acal(k, X, R)$, namely
\emph{the associated Yang-Baxter algebra}, in
the sense of Manin \cite{Manin}. The algebra
 $\Acal (k, X, R)$ is generated by
$X$ and has quadratic defining relations, $\Re(R)$
determined by $R$ similarly to
\ref{defrelations}:
\begin{equation}
\label{defrelations2}
\Re(R) =\{ (x_jx_i-c_{ij}x_{i^{\prime}}x_{j^{\prime}})  \mid
R(x_j\otimes x_i) = c_{ij} x_{i^{\prime}}\otimes
x_{j^{\prime}} \}
\end{equation}

Sometimes it will be more convenient
to work with the free associative
algebra $k\langle X\rangle,$ instead of working with the
tensor algebra, generated by $V.$
Similarly to the identification
of $X \times X$ and the set of $X^2,$
now we identify the vector spaces $V^{\otimes m}$
and $Span_kX^m$, $m \geq 1$.
 We will show that the square-free
binomial solutions of the classical Yang-Baxter equation are
closely related with a class of quadratic PBW- algebras,
the
so called \emph{skew-polynomial rings with binomial relations} and
will prove an  analogue of Theorem \ref{theoremA}. We recall the
definition.
\end{definition}
\begin{definition} \cite{T1}.
\label{skew} Let $\Acal _0=\Acal _0(k, X, \Re _0)= k<X>/ (\Re _0)$  be a finitely
presented quadratic algebra.

a) We say that $\Acal _0(k, X, \Re _0)$ is
\emph{an algebra with binomial relations of  skew-polynomial type},
if
the set of generators
$X$ is ordered: $X =\{x_1 < x_2 < \cdots < x_n\},$
and
the set of defining relations
\[
\Re _0 =\{ x_jx_i=c_{ij}x_{i^{\prime}}x_{j^{\prime}} \mid 1 \leq i
< j \leq n,\},
\]
contains precisely $n(n-1)/2$ quadratic square-free binomial
relations such that the following three conditions hold:

1) each monomial $xy$, with $x\neq y$,  $x,y \in X$ occurs in
exactly one relation in $\Re _0$; a monomial of the type $xx$ does
not occur in any relation in $\Re _0$

2) $c_{ij} \neq 0$,
for all  $i,j$ with $1\leq i < j \leq n.$

3) For every pair $i,j$ with $1\leq i < j \leq n,$
there are inequalities:
$j> i^{\prime}, \ i^{\prime} < j^{\prime}.$

b) An algebra $\Acal _0=\Acal _0(k, X, \Re _0)$
with binomial relations of
skew-polynomial type is called \emph{a skew-polynomial ring with
binomial relations} if

4) $\Re _0$ is a Groebner basis of the ideal $I = (\Re _0)$ in the
free associative algebra $k<X>$ , with respect to the
degree-lexicographic ordering of the free semigroup $\langle X
\rangle$.
\end{definition}

\begin{remark}

It follows from the Diamond Lemma, cf. \cite{B}, that
 condition \ref{skew}.4) is equivalent to each of the conditions
 $4^{\prime})$ and $4^{\prime\prime})$ below.

$4^{\prime}$) The set of ordered monomials,
\[
\mathcal{N}_0=\{ x_1^{\alpha_1}x_2^{\alpha_2}\cdots
x_n^{\alpha_n}\mid \alpha_i \geq 0, 1\leq i \leq n\}
\]
is a
k-basis of $\Acal _0$, as a k-vector space.

$4^{\prime\prime}$) The monomials $ x_kx_jx_i,$ with $k > j > i$
do not give rise to new relations in $\Acal _0.$
\end{remark}
Note that given the relations $\Re _0$, condition $4^{\prime\prime}$) is
recognizable.
\begin{definition}
Let $\Acal _0= \Acal _0(k, X, \Re _0)$ be an algebra
with binomial relations of
skew-polynomial type.

Let $V$ be the $k$-vector space with a basis $x_1, \cdots , x_n.$
Consider the linear automorphism  $R=R(\Re_0)$ of $V\otimes V$
defined as follows:

a) for each pair $i,j, 1 \leq i < j \leq n,$ we set
\[
R(x_j\otimes x_i) = c_{ij} x_{i^{\prime}}\otimes x_{j^{\prime}},
1 \leq i<j \leq n,
\]
\[
R(x_{i^{\prime}}\otimes x_{j^{\prime}}) = \frac{1}{c_{ij}} x_jx_i,
1 \leq i<j \leq n,
\]
b) for each $i, 1 \leq i \leq n$
\[
R(x_i\otimes x_i) = x_i\otimes x_i.
\]
We say that $R$ is \emph{the automorphism} associated with the
relations $\Re _0$, and denote it by $R(\Re _0)$. We also define
the bijection $r=r(\Re _0)$ of $X^2$ onto itself, as
\begin{equation}
r(xx)=xx, \text{ for all}\  x\in X,
r(x_jx_i)=(x_{i^{\prime}}x_{j^{\prime}})
\end{equation}
and
\begin{equation}
r(x_{i^{\prime}}x_{j^{\prime}}) = x_jx_i, \ \text{whenever} \
x_jx_i=c_{ij}x_{i^{\prime}}x_{j^{\prime}} \in \Re_0 .
\end{equation}
\end{definition}
\begin{lemma}
\label{skew-YBE}
Assume that $\Acal _0(k, X, \Re _0)= k\langle X\rangle/(\Re _0)$
is an algebra with
binomial relations of
skew-polynomial type, and let  $R= R(\Re _0)$ be
the automorphism of $V\otimes V$ associated with the
relations $\Re _0$. Then $R$ is a solution of the
classical Yang-Baxter equation if and only if $\Re _0$
is Groebner basis.
\end{lemma}
\begin{proof}
Assume that $R=R(\Re _0)$ is a solution of the  Yang-Baxter
equation.
We will prove that $\Re _0$ is a Groebner basis.
It will be enough to show that  each monomial
$x_kx_jx_i,$ with $k > j >i,$ can be reduced by means of
reductions defined via $\Re _0$ to a unique element
of the shape $\alpha_{ijk}x_{i^{\prime}}x_{j^{\prime}}x_{k^{\prime}},$
where
$1 \leq i^{\prime}<j^{\prime}<k^{\prime}\leq n,$
and $\alpha_{ijk}$ is a uniquely determined
coefficient, $0 \neq \alpha_{ijk} \in k$.
Let
$(X, r(R))$ be the associated set-theoretic solution,
see \ref{associatedr}.
Denote $r_1=r\times id_X, r_2=id_X\times r.$
Then the group  ${}_{gr}{\langle r_1, r_2 \rangle},$
which is isomorphic
to the symmetric group $S_3,$ acts on the set $X^3.$
Consider the orbit $\Ocal_0$ of $w=x_kx_jx_i$ under this action.
It is not difficult to see that it has precisely $6$
elements. By Lemma \ref{orderingXL3}, the relations $\Re(r)$
form a Groeber basis, therefore the orbit $\Ocal_0$
contains exactly one ordered monomial,
namely some $w_0=x_{i^{\prime}}x_{j^{\prime}}x_{k^{\prime}},$
such that
$1 \leq i^{\prime}<j^{\prime}<k^{\prime}\leq n.$

Clearly, the orbit $\Ocal$ of $x_kx_jx_i$
under the action of ${}_{gr}{\langle R^{12}, R^{23} \rangle}$
 on $kX^3$  contains the same monomials of $X^3$ as $\Ocal_0,$ but, in general,
they  occur with  non-zero coefficients which might
be different from $1.$
 In particular, $\Ocal$
contains exactly one element in \emph{normal form modulo}
$\Re _0,$  namely
$\alpha x_{i^{\prime}}x_{j^{\prime}}x_{k^{\prime}}$
where $\alpha \in k$,  $\alpha \neq 0.$ It is also clear that
each sequence of reductions  (in the sense of \cite{B})
 reduces the monomial $x_kx_jx_i$ to some element of the orbit
$\Ocal$. It follows then, that the ambiguity $x_kx_jx_i, k>j>i$
is solvable, therefore $\Re _0$ is Groebner basis.

Conversely, let  $\Re _0$ be a Groebner basis.
Consider the associated linear automorphism $R(\Re _0)$
and the associated bijective map
$r=r(R(\Re _0)): X^2\longrightarrow X^2.$
By  \cite{TM}, Theorem 1.4, $r$ is a
 solution of the set-theoretic Yang-Baxter equation.
 Now one can easily deduce
that $R(\Re _0)$ is a solution of the classical Yang-Baxter
equation.
\end{proof}
\begin{theorem}
\label{skewpolyrings} Let $V$ be finite-dimensional vector space
over a field $k$, with a $k$-basis $X.$ Suppose $R$ is a linear
automorphism of $V\otimes V$. Then the following conditions are
equivalent:
\begin{enumerate}
\item $(k, X, R)$ is a square-free binomial solution of the
classical Yang-Baxter equation. \item There exists an ordering of
$X$, $X = \{x_1 < x_2 \cdots < x_n\},$ such that the associated
 quadratic algebra
$\Acal= \Acal(k, X, R)= k\langle X \rangle/(\Re(R))$
is a
skew-polynomial ring with binomial relations.
\end{enumerate}
Furthermore, each of the above conditions implies
that $\Acal$ is a Yang-Baxter algebra  which satisfies
conditions \ref{general4} through \ref{general11} of
Theorem \ref{general}.
In particular, $\Acal$  is a Noetherian domain
and an Artin-Schelter regular ring of global dimension $n$.
\end{theorem}
\begin{proof}
$1)\Longrightarrow 2)$. Assume $(k, X, R)$ is a
a square-free binomial solution of the
classical Yang-Baxter equation. Consider the associated
set-theoretic solution   $(X, r(R)).$ It follows from \ref{theoremA}
that there exists an ordering $X = \{x_1 < \cdots < x_n\}$
such that
the relations $\Re(r(R))$ are of skew-polynomial type.
Then the relations $\Re(R)$ of the Yang-Baxter algebra
$\Acal$ associated to $(k, X, R)$ are also
of skew-polynomial type. Now Lemma \ref{skew-YBE} implies that
$\Re(R)$ is a Groebner basis, therefore $\Acal(k, X, R)$ is
a skew-polynomial ring.
The implication $1)\Longrightarrow 2)$ follows from
 Lemma \ref{skew-YBE}.

The remaining part of the theorem presents
properties
of the skew-polynomial rings with binomial relations,
$\Acal_0$, which can be extracted from our previous works.
The Noetherian properties were proved in \cite{T2},
a combinatorial proof
of the Artin-Schelter regularity of $\Acal_0$ was first
given in \cite{T3}. Conditions
\ref{general4}, through  \ref{general11} of
Theorem \ref{general},  have been deduced in \cite{TM}
from algebraic and homological properties of the semigroups
$S$ of $I$-type and the associated semigroup algebras
$kS.$
\end{proof}

{\bf Acknowledgments}. Part of this paper was written during my
visit at Harvard in 2002. It is my pleasant duty to thank the
Department of Mathematics of Harvard University for the
invitation. I express my gratitude to  David Kazhdan for
our valuable and stimulating discussions.
I thank Shahn Majid for his kind attention  to my work and
for our inspiring and fruitful  discussions.
My cordial thanks to  Michel Van Den Bergh for our
productive cooperation, for drawing my attention
to the study of set-theoretic solutions of the Yang-Baxter
equation, and for his moral support through the years.

\ifx\undefined\bysame
\newcommand{\bysame}{\leavevmode\hbox to3em{\hrulefill}\,}
\fi

\end{document}